\documentclass[11pt]{article}
\def\?{?\vadjust{\vbox to 0pt{\vss\hbox{\kern\hsize\kern1em\large\bf ?!}}}}
\usepackage[T2A]{fontenc}
\usepackage[english]{babel}
\usepackage[tbtags]{amsmath}
\usepackage{amsfonts,amssymb}
\usepackage{mathrsfs}
\usepackage{bm}
\usepackage{esint}
\usepackage{amsthm}

\newcounter{Th}[section]
\newcounter{Lm}[section]
\newcounter{Ca}[section]
\newcounter{ThA}
\newcounter{LmA}
\newcounter{Problem}[section]
\newcounter{Remark}[section]
\newcounter{Example}[section]
\newcounter{Def}[section]
\newcounter{Assum}[section]
\newcounter{Prop}[section]

\def\theTh{\arabic{section}.\arabic{Th}}
\def\theThA{\Alph{ThA}}
\def\theLmA{\Alph{LmA}}
\def\theLm{\arabic{section}.\arabic{Lm}}
\def\theCa{\arabic{section}.\arabic{Ca}}
\def\theProp{\arabic{section}.\arabic{Prop}}

\def\theRemark{\arabic{section}.\arabic{Remark}}
\def\theExample{\arabic{section}.\arabic{Example}}
\def\theDef{\arabic{section}.\arabic{Def}}

\newenvironment{Th}[1][\relax]
    {\medspace\refstepcounter{Th}{\bf Theorem \theTh.}\ \it}
    {\rm\medspace}

\newenvironment{Prop}[1][\relax]
    {\medspace\refstepcounter{Prop}{\bf Proposition \theProp.}\ \it}
    {\rm\medspace}

\newenvironment{ThA}[1][\relax]
    {\medspace\refstepcounter{ThA}{\bf Theorem \theThA.}\ \it}
    {\rm\medspace}

\newenvironment{Lm}[1][\relax]
    {\medspace\refstepcounter{Lm}{\bf Lemma \theLm.}\ \it}
    {\rm\medspace}

\newenvironment{Remark}[1][\relax]
    {\medspace\refstepcounter{Remark}{\bf Remark \theRemark.}\rm\ }
    {\medspace}

\newenvironment{Example}[1][\relax]
    {\medspace\refstepcounter{Example}{\bf Example \theExample.}\rm\ }
    {\medspace}

\newenvironment{Def}[1][\relax]
    {\medspace\refstepcounter{Def}{\bf Definition \theDef.}\rm\ }
    {\medspace}

\makeatother

\numberwithin{equation}{section}

\date{}

\begin{document}
\author{\bfseries\large A.~I.~Tyulenev
\thanks{Steklov Mathematical Institute of Russian Academy of Sciences (Moscow). E-mails:
tyulenev-math@yandex.ru, tyulenev@mi.ras.ru }}

\title{Some porosity-type properties of sets related to the $d$-Hausdorff content
}

\maketitle

\begin{abstract}
Let $S \subset \mathbb{R}^{n}$ be a nonempty set. Given $d \in [0,n)$ and a cube $\overline{Q} \subset \mathbb{R}^{n}$ with $l=l(\overline{Q}) \in (0,1]$, we show that if the $d$-Hausdorff content $\mathcal{H}^{d}_{\infty}(\overline{Q} \cap S) < \overline{\lambda}l^{d}$ for some $\overline{\lambda} \in (0,1)$, then the set $\overline{Q} \setminus S$ contains a specific cavity.
More precisely, we prove existence of a pseudometric $\rho=\rho_{S,d}$ such that for each sufficiently small $\delta > 0$
the $\delta$-neighborhood $U^{\rho}_{\delta l}(S)$ of $S$ in the pseudometric $\rho$
does not contain the whole $\overline{Q}$. Moreover, we establish the existence of constants $\overline{\delta}=\overline{\delta}(n,d,\overline{\lambda})>0$
and $\underline{\gamma}=\underline{\gamma}(n,d,\overline{\lambda})>0$ such that $\mathcal{L}^{n}(\overline{Q} \setminus U^{\rho}_{\delta l}(S)) \geq \underline{\gamma} l^{n}$ for all
$\delta \in (0,\overline{\delta})$. If, in addition, the set $S$ is $d$-lower content regular, we prove existence of
a constant $\underline{\tau}=\underline{\tau}(n,d,\overline{\lambda})>0$ such that the cube $\overline{Q}$ is $\underline{\tau}$-porous.
The sharpness of the results is illustrated by several examples.
\end{abstract}

\begin{flushleft}
\end{flushleft}
\def\B{\rlap{$\overline B$}B}
{\small\leftskip=10mm\rightskip=\leftskip\noindent}

\section{Introduction}

In many areas of geometric measure theory and geometric analysis the so-called porous sets have a significant role.
The corresponding literature is so huge that we mention only the groundbreaking papers \cite{Dol,Va,Sal,KR},
the beautiful survey \cite{Sh}, and the monograph \cite{Mat}.
Roughly speaking, $S$ is porous if, for any point $x \in S$, there are cubic holes located in $\mathbb{R}^{n} \setminus S$ arbitrary close to $x$
whose diameter is comparable with the distance to the point $x$.
Under a different nomenclature, porosity was used already in 1920 by A.~Denjoy \cite{Den}.
As far as we know, E.~P.~Dolzenko was the first who used the term ``porous set'' \cite{Dol}.

There are intimate connections between porosity properties of sets and their dimensions. It was proved in \cite{Luu}
that a set $S \subset \mathbb{R}^{n}$ is porous if, and only if, its Assouad dimension is strictly less than $n$.
The situation becomes more complicated in the context of the Hausdorff dimension.
The papers  \cite{Sal,KR,N} contain results claiming that a ``sufficiently strong porosity'' of a given set $S \subset \mathbb{R}^{n}$
implies existence of an appropriate upper bound on its Hausdorff dimension. However (in contrast to Assouad dimension), one can construct a
nonporous set $S \subset \mathbb{R}^{n}$ whose Hausdorff dimension is strictly less
than $n$.
The main goal of the present paper is to understand how the behavior of the Hausdorff contents of intersections of cubes
$Q \cap S$ with a given set $S \subset \mathbb{R}^{n}$ affects porosity-type properties of $S$.
Such sort of questions arose naturally in the study of trace problems for Sobolev spaces \cite{T2}, \cite{T3}.
This was a motivation for writing this paper.

In order to briefly describe  the main results and ideas of the present paper we fix some notation and introduce basic concepts.
As usual, given $n \in \mathbb{N}$, $\mathbb{R}^{n}$ denotes the linear space of all strings $x=(x_{1},...,x_{n})$
of real numbers. It will be convenient to equip this space with the uniform norm, i.e., $\mathbb{R}^{n}:=(\mathbb{R}^{n},\|\cdot\|_{\infty})$,
where $\|x\|_{\infty}:=\max\{|x_{1}|,...,|x_{n}|\}$. Furthermore, as usual, by $\mathcal{L}^{n}$ we denote the classical
Lebesgue measure on $\mathbb{R}^{n}$. In what follows, given a number $d \in [0,n]$ and a set $E \subset \mathbb{R}^{n}$, by $\mathcal{H}^{d}(E)$ and $\mathcal{H}^{d}_{\infty}(E)$
we will denote  the $d$-Hausdorff measure and the $d$-Hausdorff content of $E$ respectively (see the next section for the precise definitions).
Throughout the paper, the word "cube" will always mean \textit{a closed cube in $\mathbb{R}^{n}$ with sides
parallel to the coordinate axes.} Given $x \in \mathbb{R}^{n}$ and $l \geq 0$, we set $Q_{l}(x):=\prod_{i=1}^{n}[x_{i}-l,x_{i}+l]$.
In other words, $Q_{l}(x)$ is the ball centered in $x$ of radius $l$
in the space $(\mathbb{R}^{n},\|\cdot\|_{\infty})$. Given a cube $Q$, we will denote by $l(Q)$ its side length.
Given $c > 0$ and a cube $Q=Q_{l}(x)$, we let $cQ$ denote the dilation of $Q$ with respect to its
center by a factor of $c$, i.e., $cQ_{l}(x):= Q_{cl}(x)$.

We recall a slightly modified modern definition of \textit{porous sets} that is
commonly used in the literature \cite{Sh}.
First of all, given a nonempty set $S \subset \mathbb{R}^{n}$ and a parameter $\tau \in (0,1]$,
we say that \textit{a cube $Q_{l}(x)$ is $(S,\tau)$-porous}
if there is a cube $Q_{l'}(x') \subset Q_{l}(x) \setminus S$ with $l' \geq \tau l$. A cube $Q_{l}(x)$ is said to be $S$-porous if it is
$(S,\tau)$-porous for some $\tau \in (0,1]$. \textit{The family of all $(S,\tau)$-porous
cubes} will be denoted by $\mathcal{POR}_{S}(\tau)$.
Finally, we say that a set $S$ is $\tau$-porous if $Q_{l}(x) \in \mathcal{POR}_{S}(\tau)$ for all $x \in S$ and all $l \in (0,1]$.

Recall that given $d \in [0,n]$,
a closed set $S \subset \mathbb{R}^{n}$ is said to be \textit{Ahlfors--David $d$-regular} if there are constants $\operatorname{c}_{S}^{1} > 0$ and $\operatorname{c}_{S}^{2} > 0$ such that
\begin{equation}
\label{eq1.1}
\operatorname{c}_{S}^{1}l^{d} \le \mathcal{H}^{d}(Q_{l}(x) \cap S) \le \operatorname{c}_{S}^{2}l^{d} \quad \hbox{for all} \quad x \in S \quad \hbox{and all} \quad l \in (0,1].
\end{equation}
In what follows, given $d \in [0,n]$, by $\mathcal{ADR}(d)$ we denote the class of all closed Ahlfors--David $d$-regular sets.

The starting point of our investigation is the following elementary but beautiful observation made by A.~Jonsson  \cite{Jon} (see also Proposition 9.18 in \cite{Tr}).
Note that this result was an important tool in \cite{Jon,Ihn}, where traces of Besov and Lizorkin--Triebel spaces
to Ahlfors--David $d$-regular sets were studied.

\begin{ThA}
\label{ThA}
Let $d \in [0,n)$ and $S \in \mathcal{ADR}(d)$. Then there exists a constant  $\tau \in (0,\frac{1}{2})$ depending only on $d$, $n$ and $\operatorname{c}_{S}^{1}$, $\operatorname{c}_{S}^{2}$
such that $S$ is $\tau$-porous.
\end{ThA}

We should make several remarks concerning Theorem \ref{ThA}.

{\rm(\textbf{R}1)} Requirement $d < n$ is essential. Indeed, in the case $d=n$ it is obvious that $S=\mathbb{R}^{n}$ belongs to the class $\mathcal{ADR}(n)$
but the set $S$ fails to satisfy any porosity-type properties;

{\rm(\textbf{R}2)} Example 6.1 shows that an analog of Theorem \ref{ThA} fails for sets satisfying only the right-hand inequality in \eqref{eq1.1}.
One can easily show that an analog of Theorem \ref{ThA} also fails for sets satisfying only the left-hand inequality in \eqref{eq1.1}.

{\rm(\textbf{R}3)} The Ahlfors--David $d$-regularity is only a sufficient condition for the porosity of $S$ but it is
far from being necessary.

{\rm(\textbf{R}4)} Theorem \ref{ThA} has an essential drawback. Indeed, sets composed of pieces of different dimensions
do not fall into the scope of the theorem.

Recent investigations related to trace problems for Sobolev-type spaces \cite{Ry,Shv,TV,TV2,T,T2}
called for the study porosity-type properties of more complicated (in comparison with Ahlfors--David regular sets) sets that can be composed of pieces of different dimensions.
This gives a motivation for finding some less restrictive conditions on a given set $S$ which are sufficient for the existence of some porosity-type properties of $S$.

Given a number $d \in [0,n]$, a set $S \subset \mathbb{R}^{n}$ is said to be \textit{$d$-thick} or equivalently \textit{$d$-lower content regular} 
if there exists a constant $\lambda \in (0,1]$ such that
\begin{equation}
\label{eq1.1'''}
\mathcal{H}^{d}_{\infty}(Q_{l}(x) \cap S) \geq \lambda l^{d} \quad \text{for all} \quad x \in S  \quad \text{and all} \quad l \in (0,1].
\end{equation}
Since the parameter $\lambda$ will play some role below, we introduce the following notation.~Given $d \in (0,n]$ and $\lambda \in (0,1]$,
we denote by $\mathcal{LCR}(d,\lambda)$ the class of all sets $S \subset \mathbb{R}^{n}$ for which \eqref{eq1.1'''} holds.~Furthermore,
we set $\mathcal{LCR}(d):=\cup_{\lambda \in (0,1]}\mathcal{LCR}(d,\lambda)$.

As far as we know, $d$-thick sets were firstly introduced by V.~Rychkov in \cite{Ry}.
Recently, $d$-thick sets were deeply studied in \cite{AzSh}, \cite{AzVil} where those sets were called $d$-lower content regular.
The class $\mathcal{LCR}(d)$ is a natural and far reaching generalization of the class $\mathcal{ADR}(d)$.
Indeed,
\begin{equation}
\label{inclusion}
\mathcal{ADR}(d) \subset \mathcal{LCR}(d) \quad \text{for all}  \quad d \in (0,n].
\end{equation}
It is clear that $\mathcal{ADR}(n) = \mathcal{LCR}(n)$. If $d \in (0,n)$,
inclusion \eqref{inclusion} is strict. Indeed, it was noticed in \cite{Ry} and showed in \cite{TV} that any path-connected set containing at least two distinct points is $1$-thick.
On the other hand,
it is easy to built planar rectifiable curves $\Gamma \subset \mathbb{R}^{2}$ that fail to satisfy the Ahlfors-David $1$-regularity condition \cite{T} and hence, fail to satisfy the
Ahlfors-David $d$-regularity condition for any $d \in (0,2]$.
In recent papers \cite{Ry,TV,TV2,T} it was discovered that $d$-thick sets
can be effectively used in the theory of traces of function spaces.

\textbf{Problem A.} \textit{Suppose we are given parameters $d \in [0,n)$ and $\lambda \in (0,1)$. Does there exist
a constant $\tau \in (0,1/2)$ such that $\mathcal{LCR}(d,\lambda) \subset \mathcal{POR}(\tau)$?}

Unfortunately, the answer to Problem A is negative. For example, in the case $S=\mathbb{R}^{n}$
we have  $S \in \mathcal{LCR}(d,\lambda)$ for all $d \in (0,n)$ and all $\lambda \in (0,1]$ but $S \notin \mathcal{POR}(\tau)$
for any $\tau \in (0,1/2)$. The reason for that is clear. In contrast with condition \eqref{eq1.1} condition \eqref{eq1.1'''} contains
nontrivial lower bound for the corresponding content but does not contain any nontrivial upper bound.
Indeed, the trivial estimate $\mathcal{H}^{d}_{\infty}(Q_{l}(x) \cap S) \le l^{d}$ holds does not give any useful information.

\textit{Our first main result} looks like a natural generalization of Theorem \ref{ThA}.

\begin{Th}
\label{Ca.1.1}
Let $d \in [0,n)$, $\overline{\lambda} \in (0,1)$ and $\lambda \in (0,1)$.~Then there exists a constant
$\tau=\underline{\tau}(n,d,\lambda,\overline{\lambda}) \in (0,1)$ such that, for each set $S \in \mathcal{LCR}(d,\lambda)$, any cube $\overline{Q}$  with $l(\overline{Q}) \in (0,1]$ and $\mathcal{H}^{d}_{\infty}(\overline{Q} \cap S) < \overline{\lambda}(l(\overline{Q}))^{d}$ is $(S,\underline{\tau})$-porous.
\end{Th}

In fact we show that Theorem \ref{Ca.1.1} is a simple corollary of a much more deep and complicated result.~In order to formulate it we need some notation.

Recall that a pseudometric on $\mathbb{R}^{n}$
is a symmetric, nonnegative function $\rho: \mathbb{R}^{n} \times \mathbb{R}^{n} \to [0,+\infty]$
satisfying the triangle inequality. In what follows, given a pseudometric $\rho$ on $\mathbb{R}^{n}$,
we use the symbol $\mathbb{R}^{n}_{\rho}$ to denote the pseudometric space $(\mathbb{R}^{n},\rho)$. By $B^{n,\rho}_{r}(x)$
we will denote the closed ball centered at $x \in \mathbb{R}^{n}$ with radius $r$ (in the pseudometric $\rho$), i.e.,
$B^{n,\rho}_{r}(x):=\{y \in \mathbb{R}^{n}:\rho(x,y) \le r\}$.
Given a pseudometric $\rho$ on $\mathbb{R}^{n}$, a nonempty set $S \subset \mathbb{R}^{n}$ and a parameter $\tau \in (0,1]$,
we say that \textit{the ball $B^{n,\rho}_{r}(x)$ is $(S,\rho,\tau)$-porous}
if there is a ball $B^{n,\rho}_{r'}(x') \subset B^{n,\rho}_{r}(x) \setminus S$ with $r' \geq \tau r$. A ball $B$ is said to be $(S,\rho)$-porous if it is
$(S,\rho,\tau)$-porous for some $\tau \in (0,1]$.
Given $\tau \in (0,1/2)$, a set $S \subset \mathbb{R}^{n}$ is said to be $(\rho,\tau)$-porous if, for every $x \in S$ and every $r \in (0,1]$,
the ball $B^{n,\rho}_{r}(x)$ is $(S,\rho,\tau)$-porous.

Given a nonempty set $S \subset \mathbb{R}^{n}$ and parameters $d \in [0,n]$, $\lambda \in (0,1]$,
we put
\begin{equation}
\notag
\mathcal{F}_{S}(d,\lambda):=\{Q: \mathcal{H}^{d}_{\infty}(Q \cap S) \geq \lambda (l(Q))^{d}\}.
\end{equation}
Since all the cubes $Q$ are assumed to be closed, we obviously have $\{x\} \in \mathcal{F}_{S}(d,\lambda)$
for all $x \in S$.

Given parameters $d \in [0,n]$, $\lambda \in (0,1]$ and a nonempty set $S$ with $\mathcal{H}^{d}_{\infty}(S) > 0$, we define
\textit{the $(d,\lambda)$-thick with respect to $S$ distance} from an arbitrary point $y \in \mathbb{R}^{n} \setminus S$ to $S$ by the formula
\begin{equation}
\notag
\notag
\operatorname{D}_{S,d,\lambda}(y,S):=
\begin{cases}
\inf\{l(Q): Q \in \mathcal{F}_{S}(d,\lambda) \hbox{ and } y \in Q\}, \quad \text{if} \quad \{Q \in \mathcal{F}_{S}(d,\lambda):y \in Q\} \neq \emptyset;\\
+\infty, \quad \{Q \in \mathcal{F}_{S}(d,\lambda):y \in Q\}=\emptyset.
\end{cases}
\end{equation}
Given $\delta \geq 0$, we also define \textit{the $(d,\lambda)$-thick $\delta$-neighborhood} of $S$ by the formula
\begin{equation}
\label{eqq.1.2}
S_{\delta}(d,\lambda):=\{ y \in \mathbb{R}^{n}: \operatorname{D}_{S,d,\lambda}(y,S) \le \delta\}.
\end{equation}
In section 5 we introduce a pseudometric $\rho:=\rho_{S,d,\lambda}$ on $\mathbb{R}^{n}$ and show that
the set $S_{\delta}(d,\lambda)$ is a $\delta$-neighborhood of the set $S$ in the pseudometric $\rho$.

If we neglect the concrete form of holes in a given cube, one can obtain a natural
generalisation of the concept of porous cubes. Given a nonempty set $S \subset \mathbb{R}^{n}$ and a number $\gamma \in (0,1]$,
we say that a cube $Q=Q_{l}(x)$ with $x \in \mathbb{R}^{n}$ and $l > 0$
is \textit{$(S,\gamma)$-hollow} if there is a Borel set $\Omega \subset Q \setminus S$
(called \textit{the $(S,\gamma)$-cavity of $Q$}) such that $\mathcal{L}^{n}(\Omega) \geq \gamma l^{n}$.
Now we are ready to formulate \textit{our second main result}.

\begin{Th}
\label{Th1.1}
Let $d \in [0,n)$, $\overline{\lambda} \in (0,1)$ and $\lambda \in (0,1)$. Then there exist constants $\underline{\gamma}=\underline{\gamma}(n,d,\overline{\lambda}) \in (0,1]$  and $\overline{\delta}=\overline{\delta}(n,d,\lambda,\overline{\lambda}) \in (0,1)$
such that, for each set $S \subset \mathbb{R}^{n}$, for every cube $\overline{Q}$
with $l(\overline{Q}) \in (0,1]$, $\mathcal{H}^{d}_{\infty}(\overline{Q}\cap S) < \overline{\lambda}(l(\overline{Q}))^{d}$ and any $\delta \in (0,\overline{\delta}]$, the set
\begin{equation}
W_{\delta l}(\overline{Q}):=\overline{Q} \setminus S_{\delta l}(d,\lambda)
\end{equation}
is an $(S,\underline{\gamma})$-cavity of $Q$.
\end{Th}

We should make several remarks clarifying Theorem \ref{Th1.1}:

{\rm(\textbf{R}5)} It will follow from the proof that
\begin{equation}
\notag
\lim_{\overline{\lambda} \to 1}\underline{\gamma}(n,d,\overline{\lambda}) = 0,
\quad \lim_{\overline{\lambda} \to 1}\overline{\delta}(n,d,\lambda,\overline{\lambda}) = 0, \quad
 \lim_{\lambda \to 0}\overline{\delta}(n,d,\lambda,\overline{\lambda}) = 0.
\end{equation}

{\rm(\textbf{R}6)} Note that if $\overline{Q} \setminus S_{\delta l}(d,\lambda) \neq \emptyset$, then there is a ball in the pseudometric $\rho_{S,d,\lambda}$ inside $\mathbb{R}^{n} \setminus S$.
In section 5 we will show that $\rho_{S,0,\lambda}(x,y)=\|x-y\|_{\infty}$. Hence, in the case $d=0$ the condition $\overline{Q} \setminus S_{\delta l}(0,\lambda) \neq \emptyset$ is equivalent to $(S,\|\cdot\|_{\infty})$-porosity of the cube 2$\overline{Q}$.

{\rm(\textbf{R}7)} We show in Example 6.2 that condition $d < n$ is essential and cannot be dropped.

One can show that if $S$ is $(d,\lambda)$-thick,
then, for all points $x \in \mathbb{R}^{n} \setminus S$ close enough to $S$,
the $(d,\lambda)$-thick (with respect to $S$) distance
from $x$ to $S$ is comparable with the usual distance from $x$ to $S$.
This is a key observation underlying the derivation of Theorem \ref{Ca.1.1} from Theorem \ref{Th1.1}.

\textit{Our third main result} is an interesting interpretation of Theorem \ref{Th1.1}.
Recall that in different topics concerning extensions of functions from a given
closed nonempty set $S \subset \mathbb{R}^{n}$, the crucial role is played by the so-called \textit{Whitney decomposition}
of $\mathbb{R}^{n} \setminus S$ \cite{Ihn,TV,JW,Tr,T,T2,Shv}.
Recall the classical Whitney Covering Lemma \cite{St} (see Section 6.1 therein):

\textit{There exists a countable
family $W_{S}$ of closed dyadic cubes
such that:}

{\rm (i)} $\operatorname{int}Q \cap \operatorname{int}Q' = \emptyset$ \textit{for any} $Q,Q' \in W_{S}$ \textit{such that} $Q \neq Q'$;

{\rm (ii)} $\mathbb{R}^{n} \setminus S = \cup \{Q:Q \in W_{S}\}$;

{\rm (iii)} $l(Q) \le \operatorname{dist}(Q,S) \le 4l(Q)$ \textit{for all} $Q \in W_{S}$.

The family $W_{S}$ is called \textit{a Whitney decomposition of} $\mathbb{R}^{n} \setminus S$ and
cubes $Q \in W_{S}$ are called \textit{Whitney cubes.}
These cubes are closely and naturally related to the family of $S$-porous cubes.
Indeed, given a Whitney cube $Q=Q_{l_{\alpha}}(x_{\alpha}) \in W_{S}$, one can find a cube $\widetilde{Q}:=Q_{l_{\alpha}}(\widetilde{x}_{\alpha})$
whose center $\widetilde{x}_{\alpha}$ is a metric projection of $x_{\alpha}$ to $S$ such that  $Q \subset 10\widetilde{Q}$.
This proves that the cube $10\widetilde{Q}$ is $(S,\frac{1}{10})$-porous. Conversely, given an $(S,\tau)$-porous (for some $\tau \in (0,1)$)
cube $\widetilde{Q}$, one
can find a Whitney cube $Q$ such that $Q \subset c\widetilde{Q}$ for some universal constant $c \geq 1$
and, furthermore, $l(Q) \approx l(\widetilde{Q})$.

For nonexpert readers we describe informally why the Whitney cubes and porous cubes are so useful in the trace
problems for the first-order
Sobolev $W_{p}^{1}(\mathbb{R}^{n})$-spaces with $p \in (1,\infty)$. We fix a closed nonempty set $S \subset \mathbb{R}^{n}$.
In the case when $p \in (1,n]$ and $S \subset \mathbb{R}^{n}$ is regular enough or in the case
when $p > n$ and $S \subset \mathbb{R}^{n}$ is arbitrary, there are analytical tools capable of gathering
information on the behavior of a given function $f: S \to \mathbb{R}$ from any cube $Q_{l}(x)$
with $x \in S$ and $l \in (0,1]$.
Hence, one can in some sense transfer the information from porous cubes $\widetilde{Q}$ centered in $S$ to the corresponding Whitney cubes $Q$ with comparable side lengths and
then glue it smoothly using the corresponding partition of unity related to $W_{S}$. This roughly and informally explains \textit{the classical Whitney's extension method.}

Unfortunately, due to some deep analytical reasons
one cannot hope to use effectively the classical Whitney extension method for Sobolev $W_{p}^{1}(\mathbb{R}^{n})$-spaces
in the case when $1 < p \le n$ and  $S$ does not satisfy any additional regularity assumptions. For example, if
$\mathcal{H}^{d}_{\infty}(S) > 0$ for some $d \in (n-p,n]$ then only the cubes $Q \in \mathcal{F}_{S}(d,\lambda)$
(for some $\lambda \in (0,1)$) can be effectively used for gathering information about a given function $f:S \to \mathbb{R}$ \cite{TV,TV2,T2,T3}.
One cannot hope that these cubes $Q$ or even their dilated versions $cQ$ will be $S$-porous in general.
As a result, it is natural to consider an appropriate substitution for the role
of a Whitney decomposition.

Given a set $E \subset \mathbb{R}^{n}$, we define, as usual, the diameter of $E$ by letting $\operatorname{diam}E:=\sup_{x,y \in E}\|x-y\|$.
Given a family $\mathcal{G}$ of subsets of $\mathbb{R}^{n}$, by $M(\mathcal{G})$ we denote its covering multiplicity, i.e.,
the minimal $M' \in \mathbb{N}$ such that every point $x \in \mathbb{R}^{n}$ belongs to at most $M'$ sets from $\mathcal{G}$.

Now we can present the application of Theorem \ref{Th1.1}.

\begin{Th}
\label{Th.cavity.decom.}
Let $d \in (0,n)$, $c > 1$, and let $S \subset Q_{0,0}$ be a compact set with $\lambda_{S}:=\mathcal{H}^{d}_{\infty}(S) > 0$.
Then, for each
$\lambda \in (0,\frac{\lambda_{S}}{c^{d}6^{n}})$, there exist constants $\operatorname{C}_{i}=\operatorname{C}_{i}(n,\lambda,d,c)$, $i=1,2,3$, and a countable
family of Borel sets $\mathcal{W}_{S}:=\mathcal{W}_{S}(d,\lambda,c)$ such that:

{\rm (i)} $\operatorname{int}(cQ_{0,0}) \setminus S = \bigcup \{\Omega \in \mathcal{W}_{S}\}$;

{\rm (ii)} for every $\Omega \in \mathcal{W}_{S}$
\begin{equation}
\frac{1}{\operatorname{C}_{1}}\operatorname{diam}\Omega
\le \operatorname{D}_{S,d,\lambda}(\Omega,S) \le \operatorname{C}_{1}\operatorname{diam}\Omega;
\end{equation}

{\rm (iii)} for every  $\Omega \in \mathcal{W}_{S}$
\begin{equation}
\mathcal{L}^{n}(\Omega) \geq \operatorname{C}_{2}(\operatorname{diam}\Omega)^{n};
\end{equation}

{\rm (iv)} $M(\mathcal{W}_{S}) = \operatorname{C}_{3}$.

\end{Th}

\textbf{Structure of the paper.} The paper is organized as follows.

\textit{Section 2} contains an elementary background.
\textit{Section 3} which is a technical core of the paper, is based on beautiful combinatorial ideas of Yu.~Netrusov \cite{Net}.
In \textit{Section 4} we prove Theorem \ref{Ca.1.1} and Theorem \ref{Th1.1}. \textit{In section 5} given parameters $d \in (0,n)$, $\lambda \in (0,1]$
and a set $S \subset \mathbb{R}^{n}$,
we introduce a new pseudometric $\rho_{S,d,\lambda}$,
establish its basic properties and show that $(d,\lambda)$-thick $\delta$-neighborhood of the set $S \subset \mathbb{R}^{n}$
is a $\delta$-neighborhood in that pseudometric. Finally, \textit{Section 6} contains elementary examples
demonstrating the sharpness of the main theorems.

{\bf Acknowledgements.}
The author would like to thank Alexey Alimov and Roman Karasev who read first versions of this paper and made valuable remarks.
The author is grateful to the anonymous referee who find several typos.

\section{Preliminaries}

Throughout the paper $C,C_{1},C_{2},...$ will be generic positive
constants. These constants can change even
in a single string of estimates. The dependence of a constant on certain parameters is expressed, for example, by
the notation $C=C(n,p,k)$. We write $A \approx B$ if there is a constant $C \geq 1$ such that $A/C \le B \le C A$.
Given a number $c \in \mathbb{R}$ we denote by $[c]$ the integer part of $c$.

By $\mathbb{R}^{n}$ we denote the linear space of all tuples $x=(x_{1},...,x_{n})$ of real numbers \textit{equipped  with the uniform  norm} $\|\cdot\|:=\|\cdot\|_{\infty}$, i.e., $\|x\|:=\|x\|_{\infty}:=\max\{|x_{1}|,...,|x_{n}|\}$. Given a set $E \subset \mathbb{R}^{n}$, we
will denote by $\operatorname{cl}E$, $\operatorname{int}E$ and $E^{c}$
the closure, the interior, and the complement (in $\mathbb{R}^{n}$) of $E$, respectively. The symbol
$\chi_{E}$ will always mean the characteristic function of $E$. Finally, by
$\# E$ we will denote \textit{the cardinality} of $E$.
Recall that by a cube $Q \subset \mathbb{R}^{n}$ we mean a closed ball in the space $(\mathbb{R}^{n},\|\cdot\|_{\infty})$.
By \textit{a dyadic cube} we mean an arbitrary \textit{closed cube} $Q_{k,m}:=\prod_{i=1}^{n} [\frac{m_{i}}{2^{k}}, \frac{m_{i}+1}{2^{k}}]$
with $k \in \mathbb{Z}$ and $m=(m_{1},...,m_{n}) \in \mathbb{Z}^{n}$. For each $k \in \mathbb{Z}$ by $\mathcal{D}_{k}$ we denote the family of all closed dyadic cubes
with side lengths $2^{-k}$. We set
\begin{equation}
\notag
\mathcal{D}:=\cup_{k \in \mathbb{Z}}\mathcal{D}_{k}, \quad \mathcal{D}_{+}:=\cup_{k \in \mathbb{N}_{0}}\mathcal{D}_{k}.
\end{equation}
Given a family of cubes $\mathcal{Q}$ in $\mathbb{R}^{n}$ and a number $c > 1$, we set
\begin{equation}
\notag
c \mathcal{Q}:=\{c Q: Q \in \mathcal{Q}\}.
\end{equation}
Recall that, given a family $\mathcal{G}$ of subsets of $\mathbb{R}^{n}$, by $M(\mathcal{G})$ we denote its covering multiplicity, i.e.,
the minimal $M' \in \mathbb{N}$ such that every point $x \in \mathbb{R}^{n}$ belongs to at most $M'$ sets from $\mathcal{G}$.

We neeed the following elementary assertion (for details, see \cite{T2}).

\begin{Prop}
\label{Prop21}
Let $c \geq 1$ and $k \in \mathbb{N}_{0}$. Then
\begin{equation}
\label{eq.Prop21}
M(c\mathcal{D}_{k})  \le ([c]+2)^{n}.
\end{equation}
\end{Prop}


\begin{Def}
\label{Def.restriction}
Let  $\mathcal{G}=\{G_{\alpha}\}_{\alpha \in \mathcal{I}}$ be a family of sets and let $U \subset \mathbb{R}^{n}$ be a set.
We define \textit{the restriction of the family} $\mathcal{G}$ to the set $U$ by the formula
\begin{equation}
\notag
\quad \mathcal{G}|_{U}:=\{G:G \subset U\}.
\end{equation}
\end{Def}
A family $\mathcal{G}$ of sets is said to be \textit{non-overlapping} if
\begin{equation}
\notag
\operatorname{int}G \cap \operatorname{int}G' \neq \emptyset \quad \hbox{for all} \quad G,G' \in \mathcal{G} \quad \hbox{such that} \quad G \neq G'.
\end{equation}

In what follows, by \textit{a measure} we mean \textit{only a nonnegative Borel} measure on $\mathbb{R}^{n}$. By $\mathcal{L}^{n}$ we denote the classical $n$-dimensional Lebesgue measure on
$\mathbb{R}^{n}$. We say that a set $E \subset \mathbb{R}^{n}$ is \textit{measurable} if it belongs to the standard Lebesgue $\sigma$-algebra in $\mathbb{R}^{n}$.

In what follows we will commonly use  the following \textit{partial order} on the set of all
non-overlapping families of dyadic cubes.
Given two non-overlapping  families $\mathcal{Q}, \mathcal{Q'} \subset \mathcal{D}$
we write $\mathcal{Q} \succeq \mathcal{Q}'$ provided that, for every $Q' \in \mathcal{Q}'$,
there exists a unique cube $Q \in \mathcal{Q}$ such that  $Q \supset Q'$.
If, in addition, $l(Q) > l(Q')$ for all such $Q$ and $Q'$ we write
$\mathcal{Q} \succ \mathcal{Q}'$.
We say that two non-overlapping families of dyadic cubes $\mathcal{Q},\mathcal{Q}' \subset \mathcal{D}$  \textit{comparable}
if either $\mathcal{Q} \succeq \mathcal{Q}'$ or $\mathcal{Q}' \succeq \mathcal{Q}$.
Otherwise we call the corresponding families \textit{incomparable.}

Given a set $E \subset \mathbb{R}^{n}$, \textit{by a covering of the set $E$} we mean a family $\mathcal{F}$
of subsets of $\mathbb{R}^{n}$ such that $E \subset \cup\{F: F \in \mathcal{F}\}$.
Given a set $E \subset \mathbb{R}^{n}$, by \textit{a dyadic non-overlapping covering of the set} $E$ we mean a non-overlapping family $\mathcal{Q} \subset \mathcal{D}$ such that $\mathcal{Q}$ is a covering of $E$.

Given an at most countable family $\mathcal{F}$ of subsets of $\mathbb{R}^{n}$ and a number $d \in [0,n]$,
we set
\begin{equation}
\label{eqq.5.1}
\operatorname{H}^{d}(\mathcal{F}):=\sum \{(\operatorname{diam}F)^{d} : F \in \mathcal{F}\}.
\end{equation}
We also define \textit{the metric floor and the metric roof} of $\mathcal{F}$ by letting
\begin{equation}
\label{eqq.5.2}
\underline{\mu}(\mathcal{F}):=\inf \{\operatorname{diam}F : F \in \mathcal{F}\}, \quad
\overline{\mu}(\mathcal{F}):=\sup \{\operatorname{diam}F : F \in \mathcal{F}\}.
\end{equation}

In this paper we will work not only with the classical Hausdroff measures and contents but
also with their corresponding dyadic analogs.

\begin{Def}
\label{Def.content}
Let $E \subset \mathbb{R}^{n}$ be a nonempty set and $d \in [0,n]$. For any $\delta \in (0,\infty]$, we set
\begin{equation}
\label{eq2.1}
\begin{split}
&\mathcal{H}^{d}_{\delta}(E):=\inf \operatorname{H}^{d}(\mathcal{F}), \quad \mathcal{DH}^{d}_{\delta}(E):=\inf \operatorname{H}^{d}(\mathcal{Q}),
\end{split}
\end{equation}
where in the definition of $\mathcal{H}^{d}_{\delta}(E)$ the infimum is taken over all at most countable coverings $\mathcal{F}$
of the set $E$ such that $\overline{\mu}(\mathcal{F}) < \delta$ and
in the definition of $\mathcal{DH}^{d}_{\delta}(E)$ the infimum is taken over all dyadic non-overlapping coverings $\mathcal{Q}$
of the set $E$ with $\overline{\mu}(\mathcal{Q}) < \delta$.
The value $\mathcal{H}^{d}_{\infty}(E)$ is called
\textit{the $d$-Hausdorff  content} of the set $E$. The value $\mathcal{DH}^{d}_{\infty}(E)$ is called
\textit{the dyadic $d$-Hausdorff  content} of the set $E$.
We define \textit{the $d$-Hausdorff measure} and \textit{the dyadic $d$-Hausdorff measure} of the set $E$ respectively by letting
\begin{equation}
\mathcal{H}^{d}(E):=\lim_{\delta \to 0}\mathcal{H}^{d}_{\delta}(E), \quad \mathcal{DH}^{d}(E):=\lim_{\delta \to 0} \mathcal{DH}^{d}_{\delta}(E).
\end{equation}
\end{Def}

\begin{Remark}
\label{Rem2.1}
Given a set $E \subset \mathbb{R}^{n}$ and a parameter $\delta > 0$, it is easy to show that
\begin{equation}
\label{eq22}
\mathcal{H}^{d}_{\delta}(E) \le \mathcal{DH}^{d}_{\delta}(E) \le 2^{n} \mathcal{H}^{d}_{\delta}(E).
\end{equation}
\end{Remark}
\hfill$\Box$

\begin{Remark}
Let $d \in [0,n]$ and $\delta \in (0,\infty]$. Let $E \subset \mathbb{R}^{n}$ be an arbitrary set. Then by Lemma 4.6 in
\cite{Mat} and Remark \ref{Rem2.1} we get
\begin{equation}
\notag
\mathcal{H}^{d}(E)=0 \Longleftrightarrow \mathcal{DH}^{d}(E)=0
\Longleftrightarrow \mathcal{H}^{d}_{\delta}(E)=0 \Longleftrightarrow \mathcal{DH}^{d}_{\delta}(E)=0.
\end{equation}
\end{Remark}
\hfill$\Box$

\begin{Def}
\label{dyadic.cov}
Let $d \in [0,n]$ and $E \subset \mathbb{R}^{n}$ be an arbitrary set with $\mathcal{H}^{d}_{\infty}(E) > 0$. We say that a family $\mathcal{F}$
of subsets of $\mathbb{R}^{n}$ is \textit{a $d$-almost covering of the set} $E$ if there exists a set $E' \subset E$  such that $\mathcal{H}^{d}_{\infty}(E')=0$ and
$\mathcal{F}$ is a covering of $E \setminus E'$.
\end{Def}

\begin{Def}
\label{opt.cov}
Let $d \in (0,n]$ and $S \subset \mathbb{R}^{n}$ be a set with $\mathcal{H}^{d}_{\infty}(S) \in (0,+\infty)$. Given $\varepsilon > 0$,
we say that a $d$-almost covering $\mathcal{F}$ of $S$ is \textit{$\varepsilon$-optimal} if
\begin{equation}
\notag
\operatorname{H}^{d}(\mathcal{F}) \le (1+\varepsilon)\mathcal{H}^{d}_{\infty}(S).
\end{equation}
Similarly, a dyadic non-overlapping $d$-almost covering $\mathcal{Q}$ of $S$ is \textit{$\varepsilon$-optimal} if
\begin{equation}
\notag
\operatorname{H}^{d}(\mathcal{Q}) \le (1+\varepsilon) \mathcal{DH}^{d}_{\infty}(S).
\end{equation}
We say that a dyadic $\varepsilon$-optimal non-overlapping $d$-almost covering $\mathcal{Q}$
of the set $S$ is \textit{maximal}
if $\mathcal{Q} \succeq \mathcal{Q'}$ for any
$\varepsilon$-optimal dyadic non-overlapping $d$-almost covering $\mathcal{Q}'$ of $S$  \textit{comparable with} $\mathcal{Q}$.
\end{Def}

\begin{Remark}
Let $d \in (0,n]$ and $S \subset \mathbb{R}^{n}$ be a set with $\mathcal{H}^{d}_{\infty}(S) \in (0,+\infty)$. It is easy to see
that a maximal $\varepsilon$-optimal dyadic non-overlapping $d$-almost covering of the set $S$ exists, but is not unique in general.
\end{Remark}
\hfill$\Box$

\begin{Def}
\label{thick.cube}
Let $d \in [0,n]$ and $\lambda \in (0,1]$. Let $S \subset \mathbb{R}^{n}$ be a set with $\mathcal{H}^{d}_{\infty}(S) > 0$.
We say that a cube $Q=Q_{l}(x)$ with $x \in \mathbb{R}^{n}$ and $l \in [0,\infty)$ is \textit{$(d,\lambda)$-thick  with respect  to $S$} if
\begin{equation}
\label{eq3.1}
\mathcal{H}^{d}_{\infty}(Q \cap S) \geq \lambda l^{d}.
\end{equation}
Similarly,  a cube $Q \subset \mathbb{R}^{n}$ is said to be \textit{$(d,\lambda)$-dyadically thick  with respect to $S$}
if
\begin{equation}
\label{eqqq3.1}
\mathcal{DH}^{d}_{\infty}(Q \cap S) \geq \lambda l^{d}.
\end{equation}
\end{Def}

Given parameters $d \in [0,n]$, $\lambda \in (0,1]$ and a set $S \subset \mathbb{R}^{n}$, we introduce the family of all $(d,\lambda)$-thick cubes
\begin{equation}
\label{key.family}
\begin{split}
&\mathcal{F}_{S}(d,\lambda):=\{Q: \mathcal{H}^{d}_{\infty}(Q \cap S) \geq \lambda (l(Q))^{d}\}.
\end{split}
\end{equation}

\begin{Remark}
\label{Remm2.4}
Note that $Q_{0}(x)=x$ for any point $x \in \mathbb{R}^{n}$. Hence, inequality \eqref{eq3.1} trivially holds with $Q=\{x\}$ for any
$d \in (0,n]$ and $\lambda \in (0,1]$. Sometimes it will be convenient for us to consider a point $x \in S$ of a given set $S \subset \mathbb{R}^{n}$
as a $(d,\lambda)$-thick with respect to $S$ cube (whose side length is zero)
for some $d \in (0,n]$ and $\lambda \in (0,1]$.
\end{Remark}
\hfill$\Box$

In Section 3 we will currently work with a special family of dyadic cubes. This family forms some sort of building blocks for the proof
of main results of the present paper.

\begin{Def}
\label{Def4.3}
Let $S \subset \mathbb{R}^{n}$ be a nonempty set. Given $d \in [0,n]$ and $\lambda \in (0,1]$, we define \textit{the $(d,\lambda)$-keystone} for $S$ family of cubes by the formula
\begin{equation}
\notag
\mathcal{DF}_{S}(d,\lambda):=\{Q \in \mathcal{D}_{+}: \mathcal{DH}^{d}_{\infty}(Q \cap S) \geq \lambda (l(Q))^{d}\}.
\end{equation}
\end{Def}

\begin{Remark}
\label{Rem31}
Let $S \subset \mathbb{R}^{n}$ be an arbitrary set. If a cube $Q=Q_{l}(x) \in \mathcal{F}_{S}(d,\lambda)$ for some $d \in (0,n]$ and $\lambda \in (0,1]$,
then the cube $Q_{cl}(x) \in \mathcal{F}_{S}(d,\frac{\lambda}{c^{d}})$. Indeed, using the monotonicity property of the $d$-Hausdorff content we get
\begin{equation}
\notag
\mathcal{H}^{d}_{\infty}(Q_{cl}(x) \cap S) \geq \mathcal{H}^{d}_{\infty}(Q_{l}(x) \cap S) \geq \lambda l^{d}= \frac{\lambda}{c^{d}} (cl)^{d}.
\end{equation}
\end{Remark}
\hfill$\Box$

\begin{Prop}
\label{Prop.5.1}
Let $S \subset \mathbb{R}^{n}$ be a set. Let $d \in (0,n]$, $\lambda \in (0,1]$.  Then
there exists $\varepsilon_{0} > 0$ such that, for every $\varepsilon \in (0,\varepsilon_{0})$, for any cube $\overline{Q}=Q_{l}(x)$ with $\mathcal{DH}^{d}(\overline{Q} \cap S) < \lambda l^{d}$, and any $\varepsilon$-optimal dyadic nonoverlapping $d$-almost covering $\mathcal{Q}$ of the set $\overline{Q} \cap S$
\begin{equation}
\label{eq54'}
l(Q) \le 2^{k(l)}  \quad \hbox{for every} \quad Q \in \mathcal{Q},
\end{equation}
where $k(l)$ is a unique integer for which $l \in (2^{k(l)},2^{k(l)+1}]$.
\end{Prop}

\begin{proof}
We set $\overline{\lambda}:=\mathcal{DH}^{d}_{\infty}(\overline{Q} \cap S)$. By the assumptions $\overline{\lambda} < \lambda l^{d}$.
Choose $\varepsilon_{0} > 0$ so small that $(1+\varepsilon_{0})\overline{\lambda} < \lambda l^{d}$. Hence, taking $\varepsilon \in (0,\varepsilon_{0})$ and
taking an arbitrary $\varepsilon$-optimal dyadic nonoverlapping $d$-almost covering $\mathcal{Q}$ of the set $Q \cap S$, we clearly get
\begin{equation}
\notag
(l(Q))^{d} \le \operatorname{H}^{d}(\mathcal{Q}) \le (1+\varepsilon)\overline{\lambda} < \lambda l^{d} \quad \hbox{for every cube} \quad Q \in \mathcal{Q}.
\end{equation}
Since $\mathcal{Q} \subset \mathcal{D}$ we obtain \eqref{eq54'}.
\end{proof}

The following elementary proposition will be quite useful in the sequel. It exhibits relations between $(d,\lambda)$-thick cubes
and $(d,\lambda)$-dyadically thick cubes respectively.

\begin{Prop}
\label{Prop3.1}
Let $d \in (0,n)$, $\lambda \in (0,1)$ and let $S \subset \mathbb{R}^{n}$ be a Borel set with $\mathcal{H}^{d}_{\infty}(S) > 0$. Let $Q=Q_{l}(x)$ be a cube
with $l \in (0,1]$ and let $k:=[-\log_{2}l] \in \mathbb{N}_{0}$. Then,

{\rm (i)} if $cQ \in \mathcal{F}_{S}(d,\lambda)$ ($cQ \in \mathcal{DF}_{S}(d,\lambda)$) for some $c \geq 1$, then there
exists a dyadic cube $Q_{k,m} \in \mathcal{F}_{S}(d,\frac{\lambda}{([2c]+1)^{n}})$
($Q_{k,m} \in \mathcal{DF}_{S}(d,\frac{\lambda}{([2c]+1)^{n}})$) such that $Q_{k,m} \cap Q \neq \emptyset$;

{\rm (ii)} if $Q \notin \mathcal{F}_{S}(d,\frac{\lambda}{2^{d(j+1)}})$
($Q \in \mathcal{D}_{+} \setminus \mathcal{DF}_{S}(d,\frac{\lambda}{2^{d(j+1)}})$)
for some $j \in \mathbb{N}_{0}$, then every dyadic cube $Q_{k+j,m} \subset Q$ does not belong
to $\mathcal{F}_{S}(d,\lambda)$ (to $\mathcal{DF}_{S}(d,\lambda)$).

\end{Prop}

\begin{proof}
It is clear that $l \in [2^{-k},2^{-k+1})$. Since $Q$ is closed,
there are at most $([2c]+1)^{n}$ dyadic cubes $Q_{k,m}$ such that $Q_{k,m} \cap cQ \neq \emptyset$.

To prove the first claim we assume the contrary. Using the subadditivity property of $\mathcal{H}^{d}_{\infty}$ we get
\begin{equation}
\notag
\mathcal{H}^{d}_{\infty}(cQ \cap S) \le \sum\{\mathcal{H}^{d}_{\infty}(Q_{k,m} \cap S): Q_{k,m} \cap cQ \neq \emptyset\}  <   \frac{\lambda([2c]+1)^{n}}{2^{kd}([2c]+1)^{n}}  \le \lambda(l(Q))^{d}.
\end{equation}
This contradicts the assumption that $cQ \in \mathcal{F}_{S}(d,\lambda)$.

To prove the second claim assume on the contrary that there is a dyadic cube $Q_{k+j,m} \subset Q$ such that $Q_{k+j,m} \in \mathcal{F}_{S}(d,\lambda)$.
Due to monotonicity of $\mathcal{H}^{d}_{\infty}$ we get by definition of the number $k$
\begin{equation}
\notag
\mathcal{H}^{d}_{\infty}(Q \cap S) \geq \mathcal{H}^{d}_{\infty}(Q_{k+j,m} \cap S) \geq \frac{\lambda}{2^{(k+j)d}}  \geq  \frac{\lambda l^{d}}{2^{d(j+1)}}.
\end{equation}
However, this contradicts the assumption that $Q \notin \mathcal{F}_{S}(d,\frac{\lambda}{2^{d(j+1)}})$.

The corresponding dyadic analogs of the claims can be proved similarly.
\end{proof}

\section{Keystone families of cubes}

The following data are \textit{assumed to be fixed} during the whole section:

{\rm $(\operatorname{\textbf{D1}})$} arbitrary numbers $n \in \mathbb{N}$ and $d \in (0,n]$;

{\rm $(\operatorname{\textbf{D2}})$} a set $S \subset \mathbb{R}^{n}$ with $\mathcal{H}^{d}_{\infty}(S) > 0$.

Recall Definition \ref{Def4.3}. Given $\lambda \in (0,1]$, in this section we set $\mathcal{DF}(\lambda):=\mathcal{DF}_{S}(d,\lambda)$ for brevity.
In the sequel we will deal with special subfamilies of $\mathcal{DF}(\lambda)$.

\begin{Def}
\label{Def31}
Given $\lambda \in (0,1]$, we say that a family $\mathcal{Q}$ of cubes is $(d,\lambda)$-\textit{nice for the set $S$}
if the following conditions hold:

{\rm (1)} the family $\mathcal{Q}$ is a dyadic non-overlapping $d$-almost covering of $S$;

{\rm (2)} $\mathcal{Q} \subset \mathcal{DF}(\lambda)$.

\end{Def}

The following result is a modification of Lemma 2.1 of Netrusov \cite{Net} adapted to our framework.
We present a full proof to make our paper self-contained. Furthermore, we hope that the proof will clarify the driving ideas of this section.

\begin{Lm}
\label{Lm.nice.cover}
Let $\overline{Q} \in \mathcal{D}_{+}$ be such that
\begin{equation}
\label{eq4.2}
0 <\mathcal{DH}^{d}_{\infty}(\overline{Q} \cap S) < 1.
\end{equation}

Then, for each $\lambda \in (0,1)$,  there exists a family of cubes $\widehat{\mathcal{Q}}(\lambda):=\widehat{\mathcal{Q}}(\overline{Q},\lambda) \subset \mathcal{D}$ such that:

{\rm (1)} $\widehat{\mathcal{Q}}(\lambda)$ is $(d,\lambda)$-nice for $Q \cap S$;

{\rm (2)} for every cube $Q \in \widehat{\mathcal{Q}}(\lambda)$,
\begin{equation}
\label{eq32'''}
l(Q) \le \frac{l(\overline{Q})}{2};
\end{equation}

{\rm (3)} the Carleson-type packing condition,
\begin{equation}
\label{eq33'''}
\operatorname{H}^{d}(\widehat{\mathcal{Q}}(\lambda)|_{Q}) \le
(l(Q))^{d}
\end{equation}
holds for every dyadic cube $Q \subset \overline{Q}$.
\end{Lm}

\begin{proof}
Given $\lambda \in (0,1)$, we fix  $\varepsilon > 0$ so small that
\begin{equation}
\label{eq.epsilon}
0 < \tau:=  \frac{\varepsilon}{1-\lambda} < 1.
\end{equation}
We split the proof into several steps.

\textit{Step 1.}  Recall Definition \ref{opt.cov}. Given a dyadic cube $K \subset \overline{Q}$ with
\begin{equation}
\label{eqq.31}
0 < \mathcal{H}^{d}_{\infty}(K \cap S) < (l(K))^{d},
\end{equation}
let $\mathcal{Q}(K)$ be a \textit{maximal $\varepsilon$-optimal dyadic non-overlapping $d$-almost
covering of the set $K \cap S$}.
By \eqref{eqq.31} and Proposition \ref{Prop.5.1} (decreasing $\varepsilon > 0$ if necessary) we have
\begin{equation}
\label{eq46'}
l(Q) \le \frac{l(K)}{2} \quad \hbox{for every} \quad Q \in \mathcal{Q}(K).
\end{equation}
The key property of the family $\mathcal{Q}(K)$ is that the Carleson-type packing condition holds true.
More precisely, by the construction and Definition \ref{opt.cov} we have,
for every dyadic cube $Q \subset K$,
\begin{equation}
\label{eq2.6}
\operatorname{H}^{d}(\mathcal{Q}(K)|_{Q}) \le (l(Q))^{d}.
\end{equation}
Indeed, otherwise if, for some dyadic cube $Q \subset K$, inequality \eqref{eq2.6} fails, then
we modify the family $\mathcal{Q}(K)$ taking $Q$ and excluding all cubes $Q' \subset Q$, $Q' \in \mathcal{Q}(K)$. This gives an $\varepsilon$-optimal
dyadic non-overlapping $d$-almost covering of $K \cap S$. But this contradicts the maximality of $\mathcal{Q}(K)$.

\textit{Step 2.} Given a dyadic cube $K \subset \overline{Q}$, we set
\begin{equation}
\notag
\mathcal{Q}^{1}(K):=\mathcal{Q}(K) \cap \mathcal{DF}(\lambda), \quad \widetilde{\mathcal{Q}}^{1}(K):=\mathcal{Q}(K) \setminus \mathcal{Q}^{1}(K).
\end{equation}
Hence, using Definition \ref{opt.cov} and the subadditivity of $\mathcal{DH}^{d}_{\infty}$ we clearly have
\begin{equation}
\notag
\begin{split}
&\frac{1}{1+\varepsilon}\operatorname{H}^{d}(\mathcal{Q}(K))=\frac{1}{1+\varepsilon}\operatorname{H}^{d}(\mathcal{Q}^{1}(K))  +
\frac{1}{1+\varepsilon} \operatorname{H}^{d}(\widetilde{\mathcal{Q}}^{1}(K))
\le \mathcal{DH}^{d}_{\infty}(K \cap S) \\
&\le
\sum\{\mathcal{DH}^{d}_{\infty}(Q \cap S):Q \in \mathcal{Q}^{1}(K)\}+
\sum\{\mathcal{DH}^{d}_{\infty}(Q \cap S):Q \in \widetilde{\mathcal{Q}}^{1}(K)\}\\
& \le
\operatorname{H}^{d}(\mathcal{Q}^{1}(K))+\lambda\operatorname{H}^{d}(\widetilde{\mathcal{Q}}^{1}(K)).
\end{split}
\end{equation}
This clearly gives
\begin{equation}
\begin{split}
\label{eq1.2}
&\varepsilon \operatorname{H}^{d}(\mathcal{Q}^{1}(K)) \geq \Bigl(1-\lambda(1+\varepsilon)\Bigr)
\operatorname{H}^{d}(\widetilde{\mathcal{Q}}^{1}(K)).
\end{split}
\end{equation}
Hence, using \eqref{eq1.2}, Definition \ref{opt.cov} and \eqref{eq.epsilon} we get
\begin{equation}
\begin{split}
\label{eq3.4''}
&\operatorname{H}^{d}(\widetilde{\mathcal{Q}}^{1}(K)) \le \frac{\varepsilon}{(1-\lambda)(1+\varepsilon)}
\operatorname{H}^{d}(\mathcal{Q}(K)) \le \tau (l(K))^{d}.
\end{split}
\end{equation}

\textit{Step 3.}
Suppose that  we have already built,
for some $k_{0} \in \mathbb{N}$ and for every $j \in \{1,...,k_{0}\}$, families of cubes
$\mathcal{Q}^{j}=\mathcal{Q}^{j}(\overline{Q})$ and $\widetilde{\mathcal{Q}}^{j}=\widetilde{\mathcal{Q}}^{j}(\overline{Q})$ such that:

{\rm (i)} $\mathcal{Q}^{1} \subset ... \subset \mathcal{Q}^{k_{0}}$ and $\widetilde{\mathcal{Q}}^{1}\succ ... \succ\widetilde{\mathcal{Q}}^{k_{0}}$;

{\rm (ii)} $\mathcal{Q}^{k_{0}}\subset \mathcal{DF}(\lambda)|_{\overline{Q}}$;

{\rm (iii)} $\widetilde{\mathcal{Q}}^{k_{0}} \subset \mathcal{D}_{+}|_{\overline{Q}} \setminus \mathcal{DF}(\lambda)|_{\overline{Q}}$;

{\rm (iv)} the following inequality
\begin{equation}
\label{eq2.7}
\operatorname{H}^{d}(\mathcal{Q}^{k_{0}}|_{Q} \cup \widetilde{\mathcal{Q}}^{k_{0}}|_{Q}) \le (l(Q))^{d}
\end{equation}
holds for every dyadic cube $Q \subset \overline{Q}$;

{\rm (v)} it holds
\begin{equation}
\label{eqq.311}
\operatorname{H}^{d}(\widetilde{\mathcal{Q}}^{k_{0}}) \le \tau^{k_{0}}(l(\overline{Q}))^{d}.
\end{equation}

We recall notation and constructions of steps 1 and 2. We put
\begin{equation}
\begin{split}
\label{eqq3.11}
&\mathcal{Q}^{k_{0}+1}:=\cup \{\mathcal{Q}^{1}(Q):Q \in \widetilde{\mathcal{Q}}^{k_{0}}\} \cup \mathcal{Q}^{k_{0}} \quad \hbox{and}\\
&\widetilde{\mathcal{Q}}^{k_{0}+1}:=\cup \{\widetilde{\mathcal{Q}}^{1}(Q):Q \in \widetilde{\mathcal{Q}}^{k_{0}}\}.
\end{split}
\end{equation}
It is clear that conditions (i)--(iii) are satisfied with $k_{0}$ replaced by $k_{0}+1$. It remains to verify that
\eqref{eq2.7} and \eqref{eqq.311} hold with $k_{0}+1$ instead of $k_{0}$.
Indeed, an application of \eqref{eq2.6} with $K$ replaced by $Q'$ gives, for any $Q \subset \overline{Q}$,
\begin{equation}
\label{eqq.312}
\begin{split}
&\sum \{(l(Q''))^{d}: Q'' \in \mathcal{Q}^{1}(Q') \cup \widetilde{\mathcal{Q}}^{1}(Q') \hbox{ for some }
Q' \in \widetilde{\mathcal{Q}}^{k_{0}}|_{Q}\}\\
&\le \sum \{(l(Q'))^{d}: Q' \in \widetilde{\mathcal{Q}}^{k_{0}}|_{Q}\} = \operatorname{H}^{d}(\widetilde{\mathcal{Q}}^{k_{0}}|_{Q}).
\end{split}
\end{equation}
By the construction it is clear that $\widetilde{Q}^{k_{0}} \cap \mathcal{Q}^{k_{0}} = \emptyset$. Hence, combining \eqref{eq2.7}, \eqref{eqq3.11}
and \eqref{eqq.312} we get
\begin{equation}
\begin{split}
\notag
&\operatorname{H}^{d}(\mathcal{Q}^{k_{0}+1}|_{Q} \cup \widetilde{\mathcal{Q}}^{k_{0}+1}|_{Q})\\
&=
\sum \{(l(Q''))^{d}: Q' \in \widetilde{\mathcal{Q}}^{k_{0}}|_{Q} \hbox{ and }
Q'' \in \mathcal{Q}^{1}(Q') \cup \widetilde{\mathcal{Q}}^{1}(Q')\}
+ \operatorname{H}^{d}(\mathcal{Q}^{k_{0}}|_{Q})\\
&\le \operatorname{H}^{d}(\mathcal{Q}^{k_{0}}|_{Q})+
\operatorname{H}^{d}(\widetilde{\mathcal{Q}}^{k_{0}}|_{Q})=
\operatorname{H}^{d}(\mathcal{Q}^{k_{0}}|_{Q} \cup \widetilde{\mathcal{Q}}^{k_{0}}|_{Q}) \le (l(Q))^{d} \quad \text{for any dyadic cube $Q \subset K$.}
\end{split}
\end{equation}
Hence, \eqref{eq2.7} holds with $k_{0}+1$ instead of $k_{0}$.

Combining \eqref{eq3.4''}, \eqref{eqq.311}, \eqref{eqq3.11} we obtain
\begin{equation}
\operatorname{H}^{d}(\widetilde{\mathcal{Q}}^{k_{0}+1}) =
\sum\{\operatorname{H}^{d}(\mathcal{Q}^{1}(Q)): Q \in \widetilde{\mathcal{Q}}^{k_{0}}\} \le \tau \operatorname{H}^{d}(\widetilde{\mathcal{Q}}^{k_{0}}) \le
\tau^{k_{0}+1}(l(\overline{Q}))^{d}.
\end{equation}

\textit{Step 4.}
As a result, by induction we built sequences $\{\mathcal{Q}^{k}\}_{k \in \mathbb{N}}:=\{\mathcal{Q}^{k}(\overline{Q})\}_{k \in \mathbb{N}}$
and $\{\widetilde{\mathcal{Q}}^{k}\}_{k \in \mathbb{N}}:=\{\widetilde{\mathcal{Q}}^{k}(\overline{Q})\}_{k \in \mathbb{N}}$
such that conditions (i)--(v) are satisfied for any $k \in \mathbb{N}$ instead of a fixed $k_{0}$.
We set
\begin{equation}
\label{eq4.13}
\widehat{\mathcal{Q}}(\lambda):= \bigcup_{k\in\mathbb{N}}\mathcal{Q}^{k} \subset \mathcal{DF}(\lambda)|_{\overline{Q}}.
\end{equation}
Note also that according to our construction estimate \eqref{eq3.4''} implies
\begin{equation}
\label{eq1.6}
\widetilde{\mathcal{Q}}^{k} \succ \widetilde{\mathcal{Q}}^{k+1} \quad \hbox{and} \quad
\operatorname{H}^{d}(\widetilde{\mathcal{Q}}^{k}) < \tau^{k} (l(\overline{Q}))^{d} \quad \hbox{for all $k \in \mathbb{N}$}.
\end{equation}
Furthermore,
\begin{equation}
\notag
\overline{Q} \cap S \setminus \cup \{Q: Q \in \widehat{\mathcal{Q}}(\lambda)\} \subset \widetilde{\mathcal{Q}}^{k}
\quad \hbox{for all $k \in \mathbb{N}$}.
\end{equation}
Since $\tau \in (0,1)$, this leads to
\begin{equation}
\notag
\mathcal{H}^{d}_{\infty}\Bigl((\overline{Q} \cap S) \setminus \cup\{Q:Q \in \widehat{\mathcal{Q}}(\lambda)\}\Bigr) \le \varlimsup
\limits_{k \to \infty} \operatorname{H}^{d}(\widetilde{\mathcal{Q}}^{k}) = 0.
\end{equation}
Hence, by \eqref{eq4.13} the family $\widehat{\mathcal{Q}}(\lambda)$ is a dyadic non-overlapping $(d,\lambda)$-thick $d$-almost covering of the set $S \cap \overline{Q}$. This proves assertion (1) of the lemma. By
our construction, assertion (2) follows easily from \eqref{eq46'}.

Finally, it is clear from our construction that $\mathcal{Q}^{k} \subset \mathcal{Q}^{k+1}$ for all $k \in \mathbb{N}$. Combining this
fact with \eqref{eq1.6} and using inequality \eqref{eq2.7} in which $k_{0}$ is replaced by $k \in \mathbb{N}$ we get
\begin{equation}
\begin{split}
\label{eq39'}
&\operatorname{H}^{d}(\widehat{\mathcal{Q}}(\lambda)|_{Q})=\lim\limits_{k \to \infty}
\operatorname{H}^{d}(\mathcal{Q}^{k}|_{Q} \cup \widetilde{\mathcal{Q}}^{k}|_{Q}\}) \le (l(Q))^{d}
\end{split}
\end{equation}
for every dyadic cube $Q \subset \overline{Q}$. This verifies assertion (3) of the lemma.

The proof is complete.
\end{proof}

The following concept will be crucial in what follows.

\begin{Def}
\label{thick.cov}
Given $\lambda \in (0,1]$, we say that a sequence $\{\widehat{\mathcal{Q}}^{s}(\lambda)\}_{s\in \mathbb{N}_{0}}:=\{\widehat{\mathcal{Q}}^{s}_{S}(d,\lambda)\}_{s\in \mathbb{N}_{0}}$
of families of cubes \textit{is a $(d,\lambda)$-nice sequence for $S$} if
the following conditions hold:

{\rm (1)} $\widehat{\mathcal{Q}}^{0}:=\{Q \in \mathcal{D}_{0}: \mathcal{DH}^{d}_{\infty}(S \cap Q) > 0\}$;

{\rm (2)} for every $s \in \mathbb{N}$ the family $\widehat{\mathcal{Q}}^{s}(\lambda)$ is $(d,\lambda)$-nice for $S$;

{\rm (3)} $\widehat{\mathcal{Q}}^{s+1}(\lambda) \prec \widehat{\mathcal{Q}}^{s}(\lambda)$ for all $s \in \mathbb{N}_{0}$;

{\rm (4)} for each $s \in \mathbb{N}_{0}$, each $\overline{Q} \in \widehat{\mathcal{Q}}^{s}(\lambda)$ and every dyadic cube $Q \subset \overline{Q}$,
\begin{equation}
\label{eq12}
\sum\{(l(Q'))^{d}:Q' \in \widehat{\mathcal{Q}}^{s+1}(\lambda)|_{Q}\}  \le
\begin{cases}
2^{n-d}(l(Q))^{d} \quad \hbox{if } Q=\overline{Q} \hbox{ and $Q \in \mathcal{DF}(1)$};\\
(l(Q))^{d} \quad \hbox{in other cases}.
\end{cases}
\end{equation}
\end{Def}

\begin{Th}
\label{Th.nice.cover}
Given $\lambda \in (0,1)$, there exists a $(d,\lambda)$-nice for $S$ sequence
$\{\widehat{\mathcal{Q}}^{s}(\lambda)\}_{s \in \mathbb{N}_{0}}$ of families of cubes.
\end{Th}

\begin{proof}
We split the proof into two steps.

\textit{Step 1.}
We fix an arbitrary cube $\overline{Q} \in \mathcal{D}_{+}$ such that
\begin{equation}
\notag
\mathcal{DH}^{d}_{\infty}(\overline{Q} \cap S) > 0
\end{equation}
and consider \textit{two cases}.

\textit{In the first case}
\begin{equation}
\notag
\mathcal{DH}^{d}_{\infty}(\overline{Q} \cap S) < (l(\overline{Q}))^{d}.
\end{equation}
We apply Lemma \ref{Lm.nice.cover} to the cube $\overline{Q}$ and obtain a $(d,\lambda)$-nice for $\overline{Q} \cap S$ family $\widehat{\mathcal{Q}}(\overline{Q},\lambda)$ satisfying
\eqref{eq32'''} and \eqref{eq33'''}.

\textit{In the second case}
\begin{equation}
\notag
\mathcal{DH}^{d}_{\infty}(\overline{Q} \cap S) = (l(\overline{Q}))^{d},
\end{equation}
i.e., the cube $\overline{Q} \in \mathcal{DF}(1)$.
Divide $\overline{Q}$ into $2^{n}$ congruent dyadic cubes. Let $\mathcal{K}_{\overline{Q}}$ be those of them whose intersection with $S$ has
positive $\mathcal{H}^{d}_{\infty}$-content.
We put $\mathcal{K}^{g}_{\overline{Q}}:=\mathcal{K}_{\overline{Q}} \cap \mathcal{DF}(\lambda)$, $\mathcal{K}^{b}_{\overline{Q}}:=\mathcal{K}_{\overline{Q}} \setminus
\mathcal{K}^{g}_{\overline{Q}}$. For each $\overline{Q}' \in \mathcal{K}_{\overline{Q}}^{b}$
we apply Lemma \ref{Lm.nice.cover}. This gives families $\widehat{\mathcal{Q}}(\overline{Q}',\lambda)$, $\overline{Q}' \in \mathcal{K}_{\overline{Q}}^{b}$ satisfying
conditions (1) and (2) of Lemma \ref{Lm.nice.cover} in which  $\widehat{\mathcal{Q}}(\lambda)$ are replaced by $\widehat{\mathcal{Q}}(\overline{Q}',\lambda)$. We set
\begin{equation}
\notag
\widehat{\mathcal{Q}}(\overline{Q},\lambda):=\mathcal{K}_{\overline{Q}}^{g} \bigcup \Bigl(\cup\{\widehat{\mathcal{Q}}(\overline{Q}',\lambda):\overline{Q}' \in \mathcal{K}^{b}_{\overline{Q}}\}\Bigr).
\end{equation}
It is clear by the construction that
\begin{equation}
\label{eq4.15'''}
\mathcal{H}^{d}_{\infty}(\overline{Q} \cap S \setminus \cup \{Q: Q \in \widehat{\mathcal{Q}}(\overline{Q},\lambda)\})=0.
\end{equation}
Furthermore, by the construction,
\begin{equation}
\label{eq.size}
l(Q)=\frac{l(\overline{Q})}{2} \quad \text{for all} \quad Q \in \mathcal{K}^{g}_{\overline{Q}} \cup \mathcal{K}^{b}_{\overline{Q}}.
\end{equation}
Using \eqref{eq33'''} with $\widehat{\mathcal{Q}}(\overline{Q},\lambda)$ replaced by $\widehat{\mathcal{Q}}(\overline{Q}',\lambda)$, $\overline{Q}' \in \mathcal{K}^{b}_{\overline{Q}}$, taking into account
that $\#(\mathcal{K}^{g}_{\overline{Q}} \cup \mathcal{K}^{b}_{\overline{Q}})\le 2^{n}$ and finally using \eqref{eq.size}, we obtain
\begin{equation}
\begin{split}
\label{eq4.16'''}
&\operatorname{H}^{d}(\widehat{\mathcal{Q}}(\overline{Q},\lambda)) \le \sum
\{\operatorname{H}^{d}(\widehat{\mathcal{Q}}(\overline{Q}',\lambda)):\overline{Q}' \in \mathcal{K}_{\overline{Q}}^{b}\} +
\operatorname{H}^{d}(\mathcal{K}^{g}_{\overline{Q}})\\
&\le \operatorname{H}^{d}(\mathcal{K}^{b}_{\overline{Q}})+\operatorname{H}^{d}(\mathcal{K}^{g}_{\overline{Q}}) =
 \#(\mathcal{K}^{g}_{\overline{Q}} \cup \mathcal{K}^{b}_{\overline{Q}})\Bigl(\frac{l(\overline{Q})}{2}\Bigr)^{d} \le  2^{n-d}(l(\overline{Q}))^{d}.
\end{split}
\end{equation}

\textit{Step 2.} We built the desirable sequence by induction.
Clearly, by $(\operatorname{\textbf{D2}})$ the family $\widehat{\mathcal{Q}}^{0}(\lambda)$ consisting of all dyadic cubes $Q \in \mathcal{D}_{0}$ for each of which
$\mathcal{DH}^{d}_{\infty}(Q \cap S) > 0$ is nonempty. We define
\begin{equation}
\notag
\widehat{\mathcal{Q}}^{1}(\lambda):=\bigcup\{\widehat{\mathcal{Q}}(\overline{Q},\lambda): \overline{Q} \in \widehat{\mathcal{Q}}^{0}(\lambda)\} \subset \mathcal{DF}(\lambda).
\end{equation}
Suppose that we have already built, for some $j_{0} \in \mathbb{N}$, families $\widehat{\mathcal{Q}}^{0}(\lambda),...,\widehat{\mathcal{Q}}^{j_{0}}(\lambda)$
such that conditions (1)--(4) of Definition \ref{thick.cov} are satisfied for any $s \in \{0,...,j_{0}-1\}$.
Then we define
\begin{equation}
\label{eq4.17'''}
\widehat{\mathcal{Q}}^{j_{0}+1}(\lambda):=\bigcup \{\widehat{\mathcal{Q}}(\overline{Q},\lambda):\overline{Q} \in \widehat{\mathcal{Q}}^{j_{0}}(\lambda)\} \subset \mathcal{DF}(\lambda).
\end{equation}
By \eqref{eq4.15'''}, \eqref{eq4.16'''}, \eqref{eq4.17'''} conditions (1)--(4) of Definition \ref{thick.cov}
are satisfied for any $s \in \{0,...,j_{0}\}$.

As a result, by induction we get the required sequence $\{\widehat{\mathcal{Q}}^{s}(\lambda)\}_{s \in \mathbb{N}_{0}}$.

\end{proof}

Despite the fact that the proof of the following result is quite elementary, as far as we know,  it
has  never been formulated in the literature in the present form. Given $\lambda \in (0,1)$
we get a canonical decomposition of
the family $\mathcal{DF}(\lambda)$. Informally speaking, this result can be looked upon as a natural generalization of the decomposition of
the family of all dyadic cubes $\mathcal{D}_{+}$ into subfamilies $\mathcal{D}_{k}$, $k \in \mathbb{N}_{0}$.

\begin{Th}
\label{Th.adm.sys.cover.}
For each $\lambda \in (0,1)$, there exists a unique sequence  $\{\mathcal{Q}^{s}(\lambda)\}_{s \in \mathbb{N}}$ such that:

{ \rm (1)}  $\mathcal{DF}(\lambda)=\cup_{s \in \mathbb{N}}\mathcal{Q}^{s}(\lambda)$;

{ \rm (2)} for every $s \in \mathbb{N}$ the family $\mathcal{Q}^{s}(\lambda)$ is $(d,\lambda)$-nice for $S$;

{\rm (3)} $\mathcal{Q}^{s}(\lambda) \succ \mathcal{Q}^{s+1}(\lambda)$ for every $s \in \mathbb{N}$;

{\rm (4)} if, for some $\overline{Q} \in \mathcal{Q}^{s}(\lambda)$ and some $\underline{Q} \in \mathcal{Q}^{s+1}(\lambda)$, there is a
cube $Q \in \mathcal{D}_{+}$  such that
\begin{equation}
\notag
\underline{Q} \subset Q \subset \overline{Q} \quad \hbox{and} \quad l(Q) \in (l(\underline{Q}), l(\overline{Q})),
\end{equation}
then the cube $Q$ does not belong to the family $\mathcal{DF}(\lambda)$, i.e., $\mathcal{H}^{d}_{\infty}(Q \cap S) < \lambda (l(Q))^{d}$.
\end{Th}

\begin{proof} We split the proof into several steps.

\textit{Step 1.} First of all, we fix $\lambda \in (0,1)$ and for each $Q \in \mathcal{D}_{0}$
with $\mathcal{DH}^{d}_{\infty}(Q \cap S) > 0$ we denote by
$\mathcal{Q}(Q,\lambda)$ the family of all maximal dyadic cubes $Q' \in \mathcal{DF}(\lambda)$
whose side lengths are strictly less than $l(Q)$. Then, for any
$Q \in \mathcal{D}_{0}$ with $\mathcal{DH}^{d}_{\infty}(Q \cap S) > 0$ we have:

{\rm (A)} $\mathcal{Q}(Q,\lambda) \subset \mathcal{DF}(\lambda)$;

{\rm (B)} $\{Q\} \succ \mathcal{Q}(Q,\lambda)$;

{\rm (C)} the family $\mathcal{Q}(Q,\lambda)$ is $(d,\lambda)$-nice for $Q \cap S$.

Properties $(\operatorname{A})$ and $(\operatorname{B})$ are clear by the construction. To establish $(\operatorname{C})$, we apply Theorem \ref{Th.nice.cover}  and
fix a $(d,\lambda)$-nice for $S$
sequence $\{\widehat{\mathcal{Q}}^{s}(\lambda)\}_{s \in \mathbb{N}_{0}}$. Let
$j_{0} \in \mathbb{N}_{0}$ be the first number among all $j \in \mathbb{N}_{0}$ satisfying $\{Q\} \succ \{\widehat{\mathcal{Q}}^{j}(\lambda)\}$.
It is clear by the construction that $\mathcal{Q}(\lambda) \succeq \{\widehat{\mathcal{Q}}^{j_{0}}(\lambda)\}$.
Since for each $Q \in \mathcal{D}_{0}$
with $\mathcal{DH}^{d}_{\infty}(Q \cap S) > 0$ the family $\widehat{\mathcal{Q}}^{j_{0}}(\lambda)|_{Q}$ is $(d,\lambda)$-nice for the set $S \cap Q$
we complete the proof of property $(\operatorname{C})$.

\textit{Step 2.}
We build the desirable sequence $\{\mathcal{Q}^{s}(\lambda)\}_{s \in \mathbb{N}_{0}}$ by induction.

\textit{The base of induction.}
We set
\begin{equation}
\begin{split}
&\mathcal{K}^{0}_{g}:=\mathcal{D}_{0} \cap \mathcal{DF}(\lambda), \quad
\mathcal{K}^{0}_{b}:=\{Q \in \mathcal{D}_{0}: 0 < \mathcal{DH}^{d}_{\infty}(Q \cap S) < \lambda\}.
\end{split}
\end{equation}
and define
\begin{equation}
\label{eq419'}
\mathcal{Q}^{1}(\lambda):= \cup\{\mathcal{Q}(Q,\lambda): Q \in \mathcal{K}^{0}_{b}\} \cup \mathcal{K}^{0}_{g}.
\end{equation}
It follows immediately from the construction that the family $\mathcal{Q}^{1}(\lambda)$ is $(d,\lambda)$-nice for $S$.

\textit{The induction step.}
Suppose that, for some $j_{0} \in \mathbb{N}$, we have already built the families $\mathcal{Q}^{s}(\lambda)$, $s \in \{1,...,j_{0}\}$.
We put
\begin{equation}
\label{eq420'}
\mathcal{Q}^{j_{0}+1}(\lambda):=\cup\{\mathcal{Q}(Q,\lambda):Q \in \mathcal{Q}^{j_{0}}(\lambda)\}.
\end{equation}
Hence, by induction we obtain the families $\mathcal{Q}^{s}(\lambda)$ for all $s \in \mathbb{N}$.

\textit{Step 3.} It is clear that
\begin{equation}
\notag
\mathcal{Q}^{s}(\lambda) \subset \mathcal{DF}(\lambda) \quad \hbox{and} \quad \mathcal{Q}^{s}(\lambda) \succ \mathcal{Q}^{s+1}(\lambda) \quad \hbox{for all $s \in \mathbb{N}$.}
\end{equation}
Furthermore,
for each $s \in \mathbb{N}$ the family $\mathcal{Q}^{s}(\lambda)$ is $(d,\lambda)$-nice for $S$.
This proves assertions (2) and (3) of the theorem.

Suppose now that there exist $j \in \mathbb{N}$, $\overline{Q} \in \mathcal{Q}^{j}(\lambda)$, $\underline{Q} \in \mathcal{Q}^{j+1}(\lambda)$ and a cube $Q \in \mathcal{D}$ such that
\begin{equation}
\notag
\underline{Q} \subset Q \subset \overline{Q} \quad \hbox{and} \quad l(Q) \in (l(\underline{Q}),l(\overline{Q})).
\end{equation}
Note that the cube $Q \notin \mathcal{DF}(\lambda)$, since otherwise we get a contradiction
with the maximality of $\underline{Q} \in \mathcal{Q}(\overline{Q},\lambda)$.
To complete the proof it is sufficient to note that already established assertion (4) of the theorem in combination with \eqref{eq419'}, \eqref{eq420'} gives assertion (1) of the
theorem, i.e.,
\begin{equation}
\notag
\mathcal{DF}(\lambda)=\bigcup_{s \in \mathbb{N}}\mathcal{Q}^{s}(\lambda).
\end{equation}

The proof is complete.
\end{proof}

\begin{Def}
\label{Def4.4}
Given $\lambda \in (0,1)$, the sequence $\{\mathcal{Q}^{s}(\lambda)\}_{s \in \mathbb{N}}$
will be called \textit{the canonical decomposition of the family $\mathcal{DF}(\lambda)$.}
\end{Def}

The following result will be important in proving Theorem \ref{Th.52}, however we believe that it
can be interesting in itself. It reflects some interesting combinatorial properties of the
canonical decomposition $\{\mathcal{Q}^{s}(\lambda)\}_{s \in \mathbb{N}}$ of the family $\mathcal{DF}(\lambda)$.

\begin{Th}
\label{Th.33}
Let $\lambda_{1},\lambda_{2} \in (0,1)$.
Let $\{\widehat{\mathcal{Q}}^{s}(\lambda_{1})\}_{s \in \mathbb{N}}$ be a $(d,\lambda_{1})$-nice for $S$ sequence.
Let $Q \in \mathcal{D}_{+}$  and let
\begin{equation}
\notag
j_{0}:=\min\{j \in \mathbb{N}_{0}:\{Q\} \succeq \widehat{\mathcal{Q}}^{j}(\lambda_{1})|_{Q}\}.
\end{equation}
Then
\begin{equation}
\label{eq.12}
\operatorname{H}^{d}(\mathcal{C}) \le
\begin{cases}
2^{n-d}\frac{(l(Q))^{d}}{\lambda_{2}}, \quad Q \in \mathcal{DF}(1),\\
\frac{(l(Q))^{d}}{\lambda_{2}}, \quad Q \notin \mathcal{DF}(1),
\end{cases}
\end{equation}
for any family $\mathcal{C} \subset \mathcal{DF}(\lambda_{2})$ satisfying the following conditions:

{\rm (1)} $\operatorname{int}Q \cap \operatorname{int}Q' = \emptyset$ for any $Q,Q' \in \mathcal{C}$ such that $Q \neq Q'$;

{\rm (2)} $\{Q\} \succeq \mathcal{C} \succeq  \widehat{\mathcal{Q}}^{j_{0}}(\lambda_{1})|_{Q}$.

\end{Th}

\begin{proof}
In the case $\{Q\} = \widehat{\mathcal{Q}}^{j_{0}}(\lambda_{1})|_{Q}$ we clearly get $\mathcal{C}=\{Q\}$, and hence, \eqref{eq.12}
trivially holds.

Suppose now that $\{Q\} \succ \widehat{\mathcal{Q}}^{j_{0}}(\lambda_{1})|_{Q}$. Note that,
for every $Q' \in \mathcal{C}$, the family $\widehat{\mathcal{Q}}^{j_{0}}(\lambda_{1})|_{Q'}$ is $(d,\lambda_{1})$-nice for $S \cap Q'$. Hence, taking into account the inclusion $\mathcal{C} \subset \mathcal{DF}(\lambda_{2})$
and using \eqref{eq12} we obtain
\begin{equation}
\label{eq3.10''}
\begin{split}
&\operatorname{H}^{d}(\mathcal{C}) \le \frac{1}{\lambda_{2}}\sum\{\mathcal{DH}^{d}_{\infty}(Q' \cap S): Q' \in \mathcal{C}\}
\le \frac{1}{\lambda_{2}}\sum \{\operatorname{H}^{d}(\widehat{\mathcal{Q}}^{j_{0}}(\lambda_{1})|_{Q'}): Q' \in \mathcal{C}\}\\
&\le \frac{1}{\lambda_{2}}\operatorname{H}^{d}(\widehat{\mathcal{Q}}^{j_{0}}(\lambda_{1})|_{Q}) \le \frac{c}{\lambda_{2}}(l(Q))^{d},
\end{split}
\end{equation}
where $c=1$ in the case when $Q \notin \mathcal{DF}(1)$  or $c = 2^{n-d}$ in the case
when $Q \in \mathcal{DF}(1)$.

\end{proof}


\section{Main results}

Recall that given a number $\delta > 0$ and a set $S \subset \mathbb{R}^{n}$ \textit{the $\delta$-neighborhood of $S$} is defined by the formula s
\begin{equation}
\label{eqq.41}
U_{\delta}(S):=\{y \in \mathbb{R}^{n}: \|y-x\|_{\infty} < \delta \hbox{ for some } x \in S\}.
\end{equation}

We start with the following elementary observation. We recall that the metric flor and the metric roof of a given family of sets were defined in \eqref{eqq.5.2}.

\begin{Prop}
\label{P.5.1}
Let $\overline{Q}$ be an arbitrary cube in $\mathbb{R}^{n}$. Let $\mathcal{F}$ be a family of subsets in $\mathbb{R}^{n}$ such that $\overline{\mu}:=\overline{\mu}(\mathcal{F}) < \frac{l(\overline{Q})}{2}$.
Then
\begin{equation}
\label{eqq.5.3}
\mathcal{L}^{n}(\cup \{F \in\mathcal{F}: F \cap \partial \overline{Q} \neq \emptyset\}) \le 2^{2n} (l(\overline{Q}))^{n-1} \overline{\mu}.
\end{equation}
\end{Prop}

\begin{proof}
Clearly, each $n$-dimensional cube has $2^{n}$ facets. Given $i \in \{1,...,2^{n}\}$, we denote by $\partial^{i}\overline{Q}$
the $i$th facet of the cube $\overline{Q}$. It is clear that, for any $\delta > 0$, we have
\begin{equation}
\label{eqq.5.4}
U_{\delta}(\partial \overline{Q}) \subset \bigcup_{i=1}^{2^{n}}U_{\delta}(\partial^{i}\overline{Q}).
\end{equation}
Elementary geometrical observations give
\begin{equation}
\label{eqq.5.5}
\mathcal{L}^{n}(U_{\delta}(\partial^{i}\overline{Q})) \le (2\delta+l(\overline{Q}))^{n-1}2\delta \quad \hbox{for each} \quad i \in \{1,...,2^{n}\}.
\end{equation}
Combination of \eqref{eqq.5.4} and \eqref{eqq.5.5} allows to obtain
\begin{equation}
\label{eqq.5.6}
\mathcal{L}^{n}(U_{\delta}(\partial \overline{Q})) \le \sum\limits_{i=1}^{2^{n}}\mathcal{L}^{n}(U_{\delta}(\partial^{i}\overline{Q}))
\le 2^{n+1}(2\delta+l(\overline{Q}))^{n-1}\delta \quad \hbox{for every} \quad \delta > 0.
\end{equation}
By \eqref{eqq.5.2} it is clear that if $F \cap \partial \overline{Q} \neq \emptyset$ for some $F \in \mathcal{F}$
then $F \subset U_{\overline{\mu}}(\partial Q)$. Hence,
using \eqref{eqq.5.6} with $\delta=\overline{\mu}$ and taking into account that $\overline{\mu} \le 2^{-1}l(\overline{Q})$ we obtain the desirable estimate
\begin{equation}
\begin{split}
\notag
&\mathcal{L}^{n}(\cup \{F \in\mathcal{F}: F \cap \partial Q \neq \emptyset\}) \le
\mathcal{L}^{n}(U_{\overline{\mu}}(\partial \overline{Q}))\\
&\le 2^{n+1}(2\overline{\mu}+l(\overline{Q}))^{n-1}\overline{\mu}
\le  2^{2n} (l(Q))^{n-1} \overline{\mu}.
\end{split}
\end{equation}
The proof is complete.
\end{proof}

Now we require the following auxiliary result, which can be interesting in itself.

\begin{Lm}
\label{L.5.1}
Let $d \in (0,n)$, $\overline{c} > 1$ and $r > 1$.
Then there exists a number $\overline{\delta}=\overline{\delta}(n,\overline{c},r) > 0$ such that, for any $\tau > 0$, $\delta \in (0,\overline{\delta}]$
and any at most countable family $\mathcal{F}$ of subsets of $\mathbb{R}^{n}$ satisfying properties:

{\rm (1)} $\operatorname{H}^{d}(\mathcal{F}) < +\infty$;

{\rm (2)} $\delta \tau \le \underline{\mu}(\mathcal{F}) \le  \overline{\mu}(\mathcal{F}) \le \tau < +\infty$;

the following inequality
\begin{equation}
\label{eqq.5.7}
\mathcal{L}^{n}(U_{r \delta \tau}(\operatorname{F})) \le \overline{c} (\tau)^{n-d} \operatorname{H}^{d}(\mathcal{F})
\end{equation}
holds with $\operatorname{F}:=\cup\{F: F \in \mathcal{F}\}$.
\end{Lm}

\begin{proof}
We fix $\theta > 1$ so close to $1$ and choose $k^{\ast} \in \mathbb{N}$ so large that
\begin{equation}
\label{eqq.58}
\theta^{n}+\frac{([2r]+1)^{n}}{2^{k^{\ast}(n-d)}} < \overline{c}.
\end{equation}
Now  we set
\begin{equation}
\label{eqq.512'}
\overline{\delta}:=\frac{\theta-1}{2[r]+1}2^{-k^{\ast}}.
\end{equation}

We fix $\tau > 0$, $\delta \in (0,\overline{\delta}]$ and an arbitrary family $\mathcal{F}$ of subsets of $\mathbb{R}^{n}$ satisfying
properties (1), (2). We set $\underline{\mu}:=\underline{\mu}(\mathcal{F})$ and $\overline{\mu}:=\overline{\mu}(\mathcal{F})$ for brevity.
It is clear that for every set $F \in \mathcal{F}$ there is a cube $Q(F) \supset F$ with $l(Q(F))=\operatorname{diam}F$.
Such a cube is not unique in general. We fix some choice of cubes $Q(F)$, $F \in \mathcal{F}$ and define the family
\begin{equation}
\notag
\mathcal{Q}:=\mathcal{Q}(\mathcal{F}):=\{Q: Q=Q(F) \hbox{ for some } F \in \mathcal{F}\}.
\end{equation}
By our construction it is clear that
\begin{equation}
\label{eqq.5.8}
\operatorname{H}^{d}(\mathcal{Q})=\operatorname{H}^{d}(\mathcal{F}),
\quad \underline{\mu}(\mathcal{Q})=\underline{\mu},
\quad \overline{\mu}(\mathcal{Q})=\overline{\mu}.
\end{equation}
Since $l(Q) \geq \delta \tau$ it is easy to see that, for any cube $Q \in \mathcal{Q}$, there is a constant $c(Q) \in (1,[2r]+1]$ such that
\begin{equation}
\label{eqq.5.9}
U_{r\delta \tau}(\operatorname{F}) \subset \bigcup \{U_{r\delta \tau}(Q): Q \in \mathcal{Q}\} \subset \bigcup \{c(Q)Q: Q \in \mathcal{Q}\}.
\end{equation}

Our goal is to make a smart choice of the constants $c(Q)$, $Q \in \mathcal{Q}$. For this purpose we split
the family $\mathcal{Q}$ into two disjoint subfamilies. Namely, we set
\begin{equation}
\label{eqq.5.10}
\mathcal{Q}^{1}:=\{Q \in \mathcal{Q}: l(Q) > 2^{-k^{\ast}} \tau\}, \quad \mathcal{Q}^{2}:=\mathcal{Q} \setminus \mathcal{Q}^{1}.
\end{equation}
Since $\delta \in (0,\overline{\delta}]$ by \eqref{eqq.512'}
we have $U_{r\delta\tau}(Q) \subset \theta Q$ for all $Q \in \mathcal{Q}^{1}$.
On the other hand, since $l(Q) \geq \delta \tau$ for all $Q \in \mathcal{Q}$
it is clear that $U_{r\delta\tau}(Q) \subset ([2r]+1)Q$ for all $Q \in \mathcal{Q}^{2}$.
Hence, inclusion \eqref{eqq.5.9} holds with
\begin{equation}
\label{eqq.5.11}
c(Q)=
\begin{cases}
\theta, \quad \hbox{if} \quad Q \in \mathcal{Q}^{1}\\
[2r]+1, \quad \hbox{if} \quad Q \in \mathcal{Q}^{2}.
\end{cases}
\end{equation}
By \eqref{eqq.5.9}, \eqref{eqq.5.10} and \eqref{eqq.5.11} we have
\begin{equation}
\label{eqq.5.14}
\mathcal{L}^{n}(U_{r\delta \tau}(\operatorname{F}))
\le \theta^{n}\operatorname{H}^{n}(\mathcal{Q}^{1}) + ([2r]+1)^{n}\operatorname{H}^{n}(\mathcal{Q}^{2}).
\end{equation}
Since $d \in (0,n)$ by \eqref{eqq.5.2} and \eqref{eqq.5.10} we obtain the first key estimate
\begin{equation}
\label{eqq.5.15}
\operatorname{H}^{n}(\mathcal{Q}^{2}) \le \Bigl(\frac{\tau}{2^{k^{\ast}}}\Bigr)^{n-d}\operatorname{H}^{d}(\mathcal{Q}^{2}) \le \Bigl(\frac{\tau}{2^{k^{\ast}}}\Bigr)^{n-d}\operatorname{H}^{d}(\mathcal{Q})=
\Bigl(\frac{\tau}{2^{k^{\ast}}}\Bigr)^{n-d}\operatorname{H}^{d}(\mathcal{F}).
\end{equation}
Similarly, we have the second key estimate
\begin{equation}
\label{eqq.5.16}
\operatorname{H}^{n}(\mathcal{Q}^{1}) \le \tau^{n-d}\operatorname{H}^{d}(\mathcal{Q}^{1})=\tau^{n-d}\operatorname{H}^{d}(\mathcal{F}).
\end{equation}

Combining \eqref{eqq.5.14}--\eqref{eqq.5.16} and taking into account \eqref{eqq.58} we obtain \eqref{eqq.5.7} and complete
the proof.
\end{proof}

Now we are ready to prove a relatively simple result, which however will be one of the keystones in our proof the main results of the present paper.
We believe that it can be interesting in itself.
Roughly speaking, we show that if a cube $\overline{Q}$ is not $(d,\overline{\lambda})$-thick with respect to a given set $S$, then one can find
a cube $\underline{Q} \subset \overline{Q}$ such that $\mathcal{H}^{d}(\underline{Q} \cap S)/(l(\underline{Q}))^{d}$ is much smaller than $\overline{\lambda}$ but the
side length $l(\underline{Q})$ is controlled from below in a reasonable way.

Since the proof below will be quite technical, we give for the reader convenience some informal explanations of the driving ideas.
Roughly speaking, if $d \in (0,n)$ and
a cube $\overline{Q}$ is $(d,\overline{\lambda})$-thin,
then, for sufficiently large $k=k(d,\overline{\lambda}) \in \mathbb{N}$, one can construct
a family $\mathcal{F} \subset \mathcal{D}_{k}|_{\overline{Q}}$ of cardinality $\approx 2^{kn}$ such that
\begin{equation}
\notag
\sum\{\mathcal{H}^{d}_{\infty}(Q \cap S): Q \in \mathcal{F}\} \le C(n)\mathcal{H}^{d}_{\infty}(\overline{Q} \cap S)
\end{equation}
for some universal constant $C(n) \geq 1$. Hence, if we assume that $\mathcal{H}^{d}_{\infty}(Q \cap S)$ is not small for
every cube $Q \in \mathcal{F}$ and if $k$ is large enough, then taking into account that $d < n$, by elementary cardinality arguments we
get a contradiction with the smallness of $\mathcal{H}^{d}_{\infty}(\overline{Q} \cap S)$.
In order to construct $\mathcal{F}$ one should fix a small enough $\varepsilon > 0$ and fix an $\varepsilon$-optimal covering
$\mathcal{Q}$ of $\overline{Q} \cap S$. The main technical difficulty is to split the family $\mathcal{Q}$ into a ``big part'' $\mathcal{Q}^{b}$
and a ``small part'' $\mathcal{Q}^{s}$ respectively. Cubes $Q \in \mathcal{Q}^{b}$ have relatively large side lengths but the Lebesgue measure $\mathcal{L}^{n}$ of the union
of such cubes is relatively small. Hence, one should fix $k \in \mathbb{N}$ such that $2^{-k}$ is approximately equal to the minimal side length of cubes
from $\mathcal{Q}^{b}$ and then select cubes from $\mathcal{D}_{k}|_{\overline{Q}}$ that do not meet cubes $Q \in \mathcal{Q}^{b}$.


\begin{Th}
\label{Th.51}
Let $d \in (0,n)$ and $\overline{\lambda} \in (0,1)$. For each $c > 1$, there exists a constant $\underline{\kappa}=\underline{\kappa}(\overline{\lambda},n,d,c) > 0$
such that, for each compact set $S \subset Q_{0,0}$, for any cube $\overline{Q}$ satisfying
\begin{equation}
\label{eqq.415''}
l(\overline{Q}) < 1 \quad \hbox{and} \quad \mathcal{H}^{d}_{\infty}(S \cap \overline{Q}) < \overline{\lambda}(l(\overline{Q}))^{d},
\end{equation}
there exists a cube $\underline{Q} \subset \overline{Q}$
with the following properties:

{\rm (i)} $\underline{Q} \in \mathcal{D}_{+}$ and $\mathcal{H}^{d}(\underline{Q} \cap S) < (\frac{\overline{\lambda}}{c})(l(\underline{Q}))^{d}$;

{\rm (ii)} $l(\underline{Q}) \geq \underline{\kappa} l(\overline{Q})$.

\end{Th}

\begin{proof}
We fix an arbitrary compact set $S \subset Q_{0,0}$ and a cube $\overline{Q}$ satisfying \eqref{eqq.415''}. Without loss of generality we may assume that $\mathcal{H}^{d}_{\infty}(S) > 0$
because otherwise the assertion is trivial. Since $\overline{\lambda} \in (0,1)$
we fix a sufficiently small $\varepsilon \in (0,1)$ and $\overline{c} > 1$ sufficiently close to $1$ in such a way that
\begin{equation}
\label{eqq.51}
(1+\varepsilon)\overline{\lambda} < 1 \quad \hbox{and} \quad 1-\overline{c}((1+\varepsilon)\overline{\lambda})^{\frac{n}{d}} > \frac{3}{4}(1-\overline{\lambda}^{\frac{n}{d}}).
\end{equation}
We split the proof into several steps.

\textit{Step 1.} Since $\overline{Q} \notin \mathcal{F}_{S}(d,\overline{\lambda})$,
there is an $\varepsilon$-optimal
covering $\mathcal{Q}:=\mathcal{Q}_{\varepsilon}$ of the set $\overline{Q} \cap S$ such that
\begin{equation}
\label{eqq.417''}
\operatorname{H}^{d}(\mathcal{Q}) \le (1+\varepsilon)\overline{\lambda}(l(\overline{Q}))^{d}.
\end{equation}
Hence, by \eqref{eqq.415''}--\eqref{eqq.417''} we get
\begin{equation}
\begin{split}
\label{e.51}
\overline{\mu}:=\overline{\mu}(\mathcal{Q}) \le \Bigl(\operatorname{H}^{d}(\mathcal{Q})\Bigr)^{\frac{1}{d}} \le (\overline{\lambda}(1+\varepsilon))^{\frac{1}{d}}l(\overline{Q})=:\tau < 1.
\end{split}
\end{equation}

\textit{Step 2.}
Let $\overline{\delta}=\overline{\delta}(n,\overline{c},3)>0$ be the same number as in Lemma \ref{L.5.1}.
We set
\begin{equation}
\label{e.519}
\underline{\kappa}:=\min\Bigl\{\frac{\overline{\delta}}{3},\frac{\overline{\lambda}^{\frac{n}{d}}}{2^{2n+3}},\Bigl(\frac{(1-\overline{\lambda}^{\frac{n}{d}})}{c 7^{n+1}}\Bigr)^{\frac{1}{n-d}}\Bigr\}.
\end{equation}
We split the family $\mathcal{Q}$ into two disjoint subfamilies. Namely, we define a ``subfamily of big cubes'' and a ``subfamily of small cubes'' of $\mathcal{Q}$ respectively by letting
\begin{equation}
\label{eqq.519}
\mathcal{Q}^{b}_{\underline{\kappa}}:=\{Q \in \mathcal{Q}: l(Q) \geq \underline{\kappa} \tau\}, \quad \mathcal{Q}^{s}_{\underline{\kappa}}:=\mathcal{Q} \setminus \mathcal{Q}^{b}_{\underline{\kappa}}.
\end{equation}
We define
\begin{equation}
\notag
F:=\cup\{Q:Q \in \mathcal{Q}^{b}_{\underline{\kappa}}\}.
\end{equation}
The main idea is to show that \eqref{e.519} guaranties that if
\begin{equation}
\label{e.57}
k:=[-\log_{2}( \underline{\kappa} \tau)],
\end{equation}
then there are a lot of cubes from $\mathcal{D}_{k}$ inside $\overline{Q}$ that do not meet $F$ and $\partial Q$. We set
\begin{equation}
\begin{split}
\notag
&\mathcal{F}^{1}_{k}:=\{Q \in \mathcal{D}_{k}:Q \cap F \neq \emptyset\}, \quad \mathcal{F}^{2}_{k}:=\{Q \in \mathcal{D}_{k}: Q \cap \partial \overline{Q} \neq \emptyset\},\\
&\mathcal{F}^{3}_{k}:=\{Q \in \mathcal{D}_{k}:Q \subset \overline{Q} \hbox{ and } Q \notin (\mathcal{F}^{1}_{k} \cup \mathcal{F}^{2}_{k})\}.
\end{split}
\end{equation}

\textit{Step 3.}
By \eqref{e.57} we have $2^{-k} \le 2\underline{\kappa}\tau$. Hence,
$Q \subset U_{3\underline{\kappa}\tau}(F)$ for every $Q \in \mathcal{F}^{1}_{k}$.
We apply Lemma \ref{L.5.1} with $\mathcal{F}=\mathcal{Q}^{1}_{\underline{\kappa}}$, $\delta=\underline{\kappa}$, $r=3$ and take into account \eqref{eqq.417''}. This gives
\begin{equation}
\label{e.54}
\begin{split}
&V^{1}_{k}:=\mathcal{L}^{n}(\cup\{Q:Q \in \mathcal{F}^{1}_{k}\}) \le  \mathcal{L}^{n}(U_{3\underline{\kappa}\tau}(F))\\
&\le \overline{c}(\tau)^{n-d}\operatorname{H}^{d}(\mathcal{Q}) \le \overline{c}((1+\varepsilon)\overline{\lambda})^{\frac{n}{d}}(l(\overline{Q}))^{n}.
\end{split}
\end{equation}

On the other hand, using Proposition \ref{P.5.1} and taking into account the first inequality in \eqref{eqq.51}
we obtain
\begin{equation}
\begin{split}
\label{e.53}
&V^{2}_{k}:=\mathcal{L}^{n}\Bigl(\cup \{Q: Q \in \mathcal{F}_{k}^{2}\}\Bigr)
\le 2^{2n+1}(l(\overline{Q}))^{n-1}\underline{\kappa}\tau \le 2^{2n+1}\underline{\kappa}(l(\overline{Q}))^{n}.
\end{split}
\end{equation}
Using \eqref{e.519} we continue \eqref{e.53} and get
\begin{equation}
\label{eqq.424}
V^{2}_{k} \le  \frac{(\overline{\lambda})^{\frac{n}{d}}}{4}(l(\overline{Q}))^{n}.
\end{equation}

\textit{Step 4.}
By the very definition of $\mathcal{F}^{3}_{k}$ we obviously have the following fact. If $Q \in \mathcal{F}^{3}_{k}$, then
\begin{equation}
\notag
\sum \{(l(Q'))^{d}: Q' \in \mathcal{Q}^{s}_{\underline{\kappa}} \hbox{ and } Q' \cap Q \cap S \neq \emptyset\}=
 \sum \{(l(Q'))^{d}: Q' \in \mathcal{Q} \hbox{ and } Q' \cap Q \cap S \neq \emptyset\}.
\end{equation}
We use this observation and take into account that $\mathcal{Q}$ is a covering of the set $\overline{Q} \cap S$.
Hence, by definition of the Hausdorff content
it is clear that for every cube $Q \in \mathcal{F}^{3}_{k}$ we have
\begin{equation}
\label{eq.forgotten}
\mathcal{H}^{d}_{\infty}(Q \cap S) \le \sum \{(l(Q'))^{d}: Q' \in \mathcal{Q}^{s}_{\underline{\kappa}} \hbox{ and } Q' \cap Q \cap S \neq \emptyset\}.
\end{equation}
By  \eqref{eqq.519} and \eqref{e.57} it follows that
\begin{equation}
\notag
2^{-k} \geq l(Q) \quad \hbox{for all} \quad Q \in \mathcal{Q}^{s}_{\underline{\kappa}}.
\end{equation}
Hence, each cube $Q' \in \mathcal{Q}^{s}_{\underline{\kappa}}$ meets at most $3^{n}$ cubes from the family $\mathcal{F}^{3}_{k}$.
A combination of this observation with \eqref{eq.forgotten} gives
\begin{equation}
\begin{split}
\label{e.512}
&\sum\{\mathcal{H}^{d}_{\infty}(Q \cap S):Q \in \mathcal{F}^{3}_{k}\}\\
&\le \sum\limits_{Q \in \mathcal{F}^{3}_{k}}\sum \{(l(Q'))^{d}: Q' \in \mathcal{Q}^{s}_{\underline{\kappa}} \hbox{ and } Q' \cap Q \cap S \neq \emptyset\} \le 3^{n} \operatorname{H}^{d}(\mathcal{Q}^{s}_{\underline{\kappa}}).
\end{split}
\end{equation}
On the other hand, since the family $\mathcal{F}^{3}_{k}$ consists of dyadic nonoverlapping cubes with the side length $2^{-k}$
it is clear that the number of cubes in $\mathcal{F}^{3}_{k}$ can be calculated by the formula
\begin{equation}
\notag
\#\mathcal{F}^{3}_{k}= 2^{kn}\operatorname{H}^{n}(\mathcal{F}^{3}_{k})=2^{kn}\mathcal{L}^{n}(\cup\{Q:Q \in \mathcal{F}^{3}_{k}\}).
\end{equation}
From \eqref{eqq.51}, \eqref{e.54}, \eqref{e.53} it follows that
\begin{equation}
\begin{split}
\notag
&\mathcal{L}^{n}(\cup\{Q:Q \in \mathcal{F}^{3}_{k}\}) \geq \Bigl((l(\overline{Q}))^{n}-V^{1}_{k}-V^{2}_{k}\Bigr)
= \frac{(1-\overline{\lambda}^{\frac{n}{d}})}{2}(l(\overline{Q}))^{n}.
\end{split}
\end{equation}
As a result, we obtain
\begin{equation}
\label{e.513}
\#\mathcal{F}^{3}_{k} \geq 2^{kn}\frac{(1-\overline{\lambda}^{\frac{n}{d}})}{2}(l(\overline{Q}))^{n}.
\end{equation}

\textit{Step 5.} If we assume that $\mathcal{F}^{3}_{k} \subset \mathcal{DF}_{k}(d,\frac{\overline{\lambda}}{c})$,
then a combination of \eqref{eqq.417''}, \eqref{eqq.519}, \eqref{e.57}, \eqref{e.512} and \eqref{e.513} gives
\begin{equation}
\label{eqq.512}
\begin{split}
&(1+\varepsilon)\overline{\lambda}(l(\overline{Q}))^{d} \geq \operatorname{H}^{d}(\mathcal{Q}^{s}_{\underline{k}})
\geq \frac{1}{3^{n}}\frac{\overline{\lambda}}{c}2^{-kd}\#\mathcal{F}^{3}_{k}\\
&>\frac{\overline{\lambda}}{2c}
2^{k(n-d)}\frac{(1-\overline{\lambda}^{\frac{n}{d}})}{3^{n}}(l(\overline{Q}))^{n}
\geq \frac{\overline{\lambda}(1-\overline{\lambda}^{\frac{n}{d}})}{2c3^{n}}
\Bigl(\frac{l(\overline{Q})}{2\underline{\kappa}\tau}\Bigr)^{n-d}(l(\overline{Q}))^{d}.
\end{split}
\end{equation}
Using \eqref{eqq.512} and taking into account the definition of $\tau$ given in \eqref{e.51} we get
\begin{equation}
\notag
(1+\varepsilon)((1+\varepsilon)\overline{\lambda})^{\frac{n-d}{d}}\underline{\kappa}^{n-d} \geq  \frac{(1-\overline{\lambda}^{\frac{n}{d}})}{c 2^{n-d+1}3^{n}}
\end{equation}
Hence, using the first inequality in \eqref{eqq.51}, we get (recall that $\varepsilon \in (0,1)$)
\begin{equation}
\notag
\underline{\kappa}^{n-d} > \frac{(1-\overline{\lambda}^{\frac{n}{d}})}{c 7^{n+1}}.
\end{equation}
This inequality is in contradiction with \eqref{e.519}.

The proof is complete.
\end{proof}



The following concept, which was already mentioned in the introduction, gives a natural generalisation of the concept of porous cubes.

\begin{Def}
Given a set $S \subset \mathbb{R}^{n}$, a cube $Q$, and a parameter $\gamma \in (0,1]$, we say that a set
$U \subset Q$ is an $(S,\gamma)$-cavity of the cube $Q$ if
\begin{equation}
\notag
U \subset Q \setminus S \quad \hbox{and} \quad \mathcal{L}^{n}(U \setminus S) \geq \gamma (l(Q))^{n}.
\end{equation}
We say that $Q$ is $(S,\gamma)$-hollow if there exists an $(S,\gamma)$-cavity $U$ of the cube $Q$.
\end{Def}

We need some notation. Given numbers  $d \in (0,n)$, $\lambda \in (0,1)$, $r \geq 1$
and a set $S \subset \mathbb{R}^{n}$ with $\mathcal{H}^{d}_{\infty}(S) > 0$, we define,
for each $\varkappa \in (0,1]$ and any cube $Q \subset \mathbb{R}^{n}$, the set
\begin{equation}
\begin{split}
\notag
&U_{\varkappa}(Q,r):=U_{\varkappa}(Q,d,\lambda,r):= Q \setminus \Bigl(\cup\{r Q': Q' \in \mathcal{DF}_{S}(d,\lambda), l(Q') \le \varkappa l(Q)\}\Bigr).\\
\end{split}
\end{equation}

\begin{Th}
\label{Th.52}
Let $d \in (0,n)$, $\lambda \in (0,1)$ and $r \geq 1$.
Then, for each $\gamma \in (0,1-2^{d-n})$, there exists a number
$\overline{\varkappa} =  \overline{\varkappa}(\gamma,n,d,\lambda,r) \in (0,1)$ such that, for every
Borel set $S \subset \mathbb{R}^{n}$  and every cube $\overline{Q}=Q_{k,m} \in \mathcal{D}_{+}$ with
\begin{equation}
\mathcal{DH}^{d}_{\infty}(\overline{Q} \cap S) < (l(\overline{Q}))^{d}
\end{equation}
the sets $U_{\varkappa}(\overline{Q},d,\lambda,r)$ are $(S,\gamma)$-cavities of the cube $\overline{Q}$ for all $\varkappa \in (0,\overline{\varkappa})$.
\end{Th}

\begin{proof}
We fix a Borel set $S \subset \mathbb{R}^{n}$ and a cube $\overline{Q}=Q_{k,m} \in \mathcal{D}_{+} \notin \mathcal{DF}_{S}(d,1)$.
Without loss of generality we assume that $\mathcal{H}^{d}_{\infty}(S) > 0$ because otherwise the assertion is trivial.
We also fix a parameter  $\gamma \in (0,1-2^{d-n})$.
During the proof we write for brevity $\mathcal{DF}:=\mathcal{DF}_{S}(d,\lambda)$.
Let $\{\widehat{\mathcal{Q}}^{s}\}_{s \in \mathbb{N}}$ be an arbitrary $(d,\lambda)$-nice for $S$
sequence. We split the proof into several steps.

\textit{Step 1.} We fix a number $\overline{c} > 1$ so close to $1$ that
\begin{equation}
\label{eqq.5.29}
\frac{\overline{c}}{2^{n-d}} < \frac{1}{2}(1-\gamma + 2^{d-n}).
\end{equation}
Let $\overline{\delta}=\overline{\delta}(\overline{c},n,r)$ be the same as in Lemma \ref{L.5.1}. Now we fix the minimal $k \in \mathbb{N}$ for which
\begin{equation}
\label{eqq.5.30}
2^{-k} < \overline{\delta} \quad \hbox{and} \quad \frac{r^{n}}{2^{k(n-d)}\lambda} + \frac{2^{2n}}{2^{k}} < \frac{1}{2}(1-\gamma-2^{d-n}).
\end{equation}

\textit{Step 2.}
We define
\begin{equation}
\label{eqq.5.31}
s_{0}:=\min\{s \in \mathbb{N}_{0}:\{Q\} \succ \widehat{\mathcal{Q}}^{s+1}|_{\overline{Q}}\}.
\end{equation}
Since $Q \in \mathcal{D}_{+} \notin \mathcal{DF}_{S}(d,1)$ by Definition \ref{thick.cov} we have
\begin{equation}
\label{e.518}
\operatorname{H}^{d}(\widehat{\mathcal{Q}}^{s_{0}+1}|_{\overline{Q}})  \le (l(\overline{Q}))^{n-d}.
\end{equation}

\textit{Step 3.}
We introduce the family
\begin{equation}
\notag
\mathcal{K}:=\{Q \in \mathcal{DF}: r Q \cap  \overline{Q} \neq \emptyset \hbox{ and } l(Q) \le 2^{-k-1}l(\overline{Q})\}.
\end{equation}
We split $\mathcal{K}$ into three subfamilies. More precisely, we set
\begin{equation}
\label{e.517}
\begin{split}
&\mathcal{K}^{1}:=\{Q : Q\in \mathcal{K} \hbox{ and } \operatorname{int}Q \subset \mathbb{R}^{n}\setminus \overline{Q}\}, \\
&\mathcal{K}^{2}:=\{Q \in \mathcal{K} \setminus \mathcal{K}^{1}: \exists \widehat{Q} \in \widehat{\mathcal{Q}}^{s_{0}+1} \hbox{ s.t. }
\widehat{Q} \subset Q\},\\
&\mathcal{K}^{3}:=\{Q \in \mathcal{K} \setminus \mathcal{K}^{1}: \exists \widehat{Q} \in \widehat{\mathcal{Q}}^{s_{0}+1}
\hbox{ s.t. } Q \subset \widehat{Q}
\hbox{ and } l(Q) < l(\widehat{Q})\}.\\
\end{split}
\end{equation}
It follows directly from the construction that
\begin{equation}
\label{eq.328}
\mathcal{K} \subset \mathcal{K}^{1} \cup \mathcal{K}^{2} \cup \mathcal{K}^{3}.
\end{equation}

\textit{Step 4.}
Using Proposition \ref{P.5.1}, we get
\begin{equation}
\label{eqq.5.35}
\mathcal{L}^{n}\Bigl(\cup\{(r Q \cap \overline{Q}): Q \in \mathcal{K}^{1}\}\Bigr) \le 2^{2n}(l(\overline{Q}))^{n}\frac{1}{2^{k}}.
\end{equation}

\textit{Step 5.}
Let $\overline{\mathcal{K}}^{2} \subset \mathcal{K}^{2}$ be the family of all maximal (with respect to inclusion)
dyadic cubes from the family $\mathcal{K}^{2}$.
We obviously get
\begin{equation}
\label{eqq.5.36}
\bigcup\{Q: Q \in \mathcal{K}^{2}\} \subset \bigcup\{Q: Q \in \overline{\mathcal{K}}^{2}\}.
\end{equation}
Consider also the family
\begin{equation}
\notag
\mathcal{C}:=\overline{\mathcal{K}}^{2} \cup \{Q \in \widehat{Q}^{s_{0}+1}|_{\overline{Q}}: \operatorname{int}Q \cap \operatorname{int}Q' = \emptyset \hbox{ for all } Q' \in \overline{\mathcal{K}}^{2}\}.
\end{equation}
We use \eqref{eqq.5.36}, then take into account that $l(Q) \le 2^{-k}l(\overline{Q})$ for all $Q \in \mathcal{Q}^{2}_{k}$
and finally apply Theorem \ref{Th.33}. As a result, we obtain
\begin{equation}
\label{eqq.5.37}
\begin{split}
&\mathcal{L}^{n}\Bigl(\cup\{r Q: Q \in \mathcal{K}^{2}\}\Bigr) \le \mathcal{L}^{n}
\Bigl(\cup\{r Q: Q \in \overline{\mathcal{K}}^{2}\}\Bigr) \le r^{n} \operatorname{H}^{n}(\overline{\mathcal{K}}^{2})\\
&\le
r^{n}\Bigl(\frac{l(\overline{Q})}{2^{k}}\Bigr)^{n-d}\operatorname{H}^{d}(\overline{\mathcal{K}}^{2}) \le
r^{n}\Bigl(\frac{l(\overline{Q})}{2^{k}}\Bigr)^{n-d}\operatorname{H}^{d}(\mathcal{C})
\le \frac{r^{n}}{2^{k(n-d)}\lambda}(l(\overline{Q}))^{n}.
\end{split}
\end{equation}

\textit{Step 6.}
Let $\widehat{\mathcal{Q}} \subset \widehat{\mathcal{Q}}^{s_{0}+1}|_{\overline{Q}}$ be the family consisting of all cubes $\widehat{Q}$
for each of which there exists a cube $Q \in \mathcal{K}^{3}$ such that $Q \subset \widehat{Q}$ and
$l(Q) < l(\widehat{Q})$. Letting $\tau:=2^{-1} l(\overline{Q})$ by \eqref{eqq.5.31} we get
\begin{equation}
\overline{\mu}(\widehat{\mathcal{Q}}) \le \overline{\mu}(\widehat{\mathcal{Q}}^{s_{0}+1}|_{\overline{Q}}) \le \tau.
\end{equation}
Using the first inequality in \eqref{eqq.5.30} we have
\begin{equation}
\label{eqq.5.38}
\cup\{r Q: Q \in \mathcal{K}^{3}\} \subset U_{r\overline{\delta}\tau}(\cup\{\widehat{Q}:\widehat{Q} \in \widehat{\mathcal{Q}}\}).
\end{equation}
We use \eqref{eqq.5.38}, then apply Lemma \ref{L.5.1}, and finally use \eqref{e.518}. This yields
\begin{equation}
\label{eqq.5.39}
\begin{split}
&\mathcal{L}^{n}(\cup\{r Q: Q \in \mathcal{K}^{3}\}) \le
\mathcal{L}^{n}(U_{r\overline{\delta}\tau}(\cup\{\widehat{Q}:\widehat{Q} \in \widehat{\mathcal{Q}}\}))\\
&\le \overline{c}\Bigl(\frac{l(\overline{Q})}{2}\Bigr)^{n-d}
\operatorname{H}^{d}(\widehat{\mathcal{Q}}) \le \overline{c}\Bigl(\frac{l(\overline{Q})}{2}\Bigr)^{n-d}
\operatorname{H}^{d}(\widehat{\mathcal{Q}}^{s_{0}+1}|_{\overline{Q}}) \le \frac{\overline{c}}{2^{n-d}}(l(\overline{Q}))^{n}.
\end{split}
\end{equation}

\textit{Step 7.} We set $\overline{\varkappa}:=2^{-k-1}$.
Collecting \eqref{eqq.5.29}, \eqref{eqq.5.30}, \eqref{eqq.5.35}, \eqref{eqq.5.37}, \eqref{eqq.5.39} we obtain
\begin{equation}
\label{eq3.21''}
\begin{split}
&\mathcal{L}^{n}(\overline{Q} \setminus U_{\overline{\varkappa}}(\overline{Q},d,\lambda,r)) \le (1-\gamma)(l(\overline{Q}))^{n}. \\
\end{split}
\end{equation}
Taking into account that
$U_{\overline{\varkappa}}(\overline{Q},d,\lambda,r) \subset U_{\varkappa}(\overline{Q},d,\lambda,r)$ for all $\varkappa \in (0,\overline{\varkappa})$
we complete the proof.

\end{proof}

Now we are ready to prove \textit{the second main result} of the present paper. We recall that the $(d,\lambda)$-thick $\delta$-neighborhood of a given set $S \subset \mathbb{R}^{n}$
was defined in \eqref{eqq.1.2}.

\textit{Proof of Theorem \ref{Th1.1}.}
The case $d=0$ is trivial. Hence, we assume that $d \in (0,n)$. Fix an arbitrary set $S \subset \mathbb{R}^{n}$ and an arbitrary cube $\overline{Q}=Q_{l}(x)$ satisfying assumptions of the theorem.
An application of Theorem \ref{Th.51} with $c=2^{n}$ gives an existence of a constant
$\underline{\kappa}:=\underline{\kappa}(\overline{\lambda},n,d,2^{n})$ and a
cube $\underline{Q} \in \mathcal{D}_{+}$ with
side length
\begin{equation}
\label{eqq.4.44}
l(\underline{Q}) \geq \underline{\kappa}l(\overline{Q})
\end{equation}
such that
$\mathcal{H}^{d}_{\infty}(\underline{Q} \cap S) < \frac{\overline{\lambda}}{2^{n}}$. By Remark \ref{Rem2.1},
\begin{equation}
\notag
\mathcal{DH}^{d}_{\infty}(\underline{Q}\cap S) < \overline{\lambda} l(\underline{Q}).
\end{equation}
Hence, by Theorem \ref{Th.52} there exists a constant $\overline{\varkappa}:=\overline{\varkappa}(\frac{1-2^{d-n}}{2},n,d,\frac{\lambda}{3^{n}},3)$ such that
\begin{equation}
\label{eqq.4.45}
\mathcal{L}^{n}(U_{\overline{\varkappa}}(\underline{Q},d,\frac{\lambda}{3^{n}},3)) \geq \Bigl(\frac{1-2^{d-n}}{2}\Bigr)(l(\underline{Q}))^{n}.
\end{equation}

By Proposition \ref{Prop3.1} for any cube $Q \in \mathcal{F}_{S}(d,\lambda)$ there is a cube $Q' \in \mathcal{D}_{j}$ with $j=-[\log_{2}l(Q)]$ such that
$Q' \in \mathcal{DF}_{S}(d,\lambda/3^{n})$. Clearly, $Q \subset 3 Q'$. Hence, if we set
\begin{equation}
\notag
\overline{\delta}:=\overline{\delta}(n,d,\overline{\lambda},\lambda):=\underline{\kappa}(\overline{\lambda},n,d,2^{n})
\overline{\varkappa}(\frac{1-2^{d-n}}{2},n,d,\frac{\lambda}{3^{n}},3)=\underline{\kappa}\overline{\varkappa}
\end{equation}
then for any $\delta \in (0,\overline{\delta})$
\begin{equation}
\label{eqq.incl'}
W_{\delta l}(\overline{Q},d,\lambda) \supset U_{\overline{\varkappa}}(\underline{Q},d,\frac{\lambda}{3^{n}},3).
\end{equation}

Now we set
\begin{equation}
\notag
\underline{\gamma}(\overline{\lambda},n,d) = \Bigl(\frac{1-2^{d-n}}{2}\Bigr)(\underline{\kappa}(\overline{\lambda},n,d,2^{n}))^{n}.
\end{equation}
As a result, by \eqref{eqq.4.44}, \eqref{eqq.4.45}, \eqref{eqq.incl'} we deduce
\begin{equation}
\mathcal{L}^{n}(W_{\delta l}(\overline{Q},d,\lambda)) \geq \mathcal{L}^{n}(U_{\overline{\varkappa}}(\underline{Q},d,\frac{\lambda}{3^{n}},3))
\geq  \underline{\gamma}(\overline{\lambda},n,d)(l(\overline{Q}))^{n}.
\end{equation}

The proof is complete.
\hfill$\Box$

\begin{Remark}
\label{Rem.simpleobserv}
It is easy to show that if $d \in (0,n]$, $\overline{\lambda} \in (0,1)$, $S \subset \mathbb{R}^{n}$ is a nonempty set, and a cube $\overline{Q}=Q_{l}(x)$
with $l \in (0,1]$ is such that $\overline{Q} \notin \mathcal{F}_{S}(d,\overline{\lambda})$, then the cube $\overline{Q}$ is $(S,1-\overline{\lambda}^{\frac{n}{d}})$-hollow.

Indeed, by Definition \ref{Def.content} there is an at most countable covering $\mathcal{U}$
of the set $\overline{Q} \cap S$ such that
\begin{equation}
\notag
\operatorname{H}^{d}(\mathcal{U}) < \overline{\lambda} l^{d}.
\end{equation}
For each set $U \in \mathcal{U}$ there is a cube $Q(U) \supset U$ with $l(Q)=\operatorname{diam}U$.
It is clear that
\begin{equation}
\notag
l(Q(U)) < \overline{\lambda}^{\frac{1}{d}} l \quad \hbox{for all} \quad U \in \mathcal{U}.
\end{equation}
This gives
\begin{equation}
\label{eq5.16}
\operatorname{H}^{n}(\mathcal{U}) < \overline{\lambda}^{\frac{n}{d}-1} l^{n-d}
\operatorname{H}^{d}(\mathcal{U}) \le \overline{\lambda}^{\frac{n}{d}} l^{n}.
\end{equation}
Since $\overline{\lambda} < 1$, the required result follows from \eqref{eq5.16} and the subadditivity property of the Lebesgue measure $\mathcal{L}^{n}$.

It is clear that there is a huge difference between the elementary observation given above and Theorem \ref{Th1.1}. The former observation
does not give any information about the structure of cavities in cubes whose intersection with $S$ have relatively small $d$-Hausdorff content.
On the other hand,  informally speaking Theorem \ref{Th1.1} claims that the corresponding cavities in cubes are located at some ``nonzero depth'' in $\mathbb{R}^{n} \setminus S$
with respect to the special distance.
\end{Remark}
\hfill$\Box$

Now we show that Theorem \ref{Th1.1}, which was called the first main result of the present paper admits a significant clarification in the context of $d$-thick sets.

\textit{Proof of Theroem \ref{Ca.1.1}.} The case $d=0$ is trivial. Hence, we assume that $d \in (0,n)$. Since $S$ is $(d,\lambda)$-thick we have
\begin{equation}
\notag
S_{\varepsilon}(d,\lambda) \supset U_{\varepsilon}(S) \quad \hbox{for all} \quad \varepsilon \in (0,1].
\end{equation}
By Theorem \ref{Th1.1} this implies that
\begin{equation}
\notag
\overline{Q} \setminus U_{\delta l}(S) \neq \emptyset \quad \hbox{for all} \quad \delta \in (0,\overline{\delta}].
\end{equation}
Hence, letting $\underline{\tau}:=\underline{\tau}(n,d,\overline{\lambda},\lambda)=\frac{1}{2}\overline{\delta}(n,d,\overline{\lambda},\lambda)$, and taking an arbitrary point $x^{\ast} \in \overline{Q} \setminus U_{\delta l}(S)$ we obtain
\begin{equation}
\notag
Q_{\underline{\tau}}(x^{\ast}) \subset \overline{Q} \setminus S.
\end{equation}
The proof is complete.
\hfill$\Box$

\section{Applications}

In this section we present the proof of Theorem \ref{Th.cavity.decom.}. Furthermore, we introduce some new concepts, which can be of independent
interest.

The following data are assumed to be fixed during the whole section:

{\rm (\textbf{A})}  arbitrary numbers $n \in \mathbb{N}$ and $d \in [0,n)$;

{\rm (\textbf{B})} a compact set $S \subset Q_{0,0}$ with $\lambda_{S}:=\mathcal{H}^{d}_{\infty}(S) > 0$.

Recall \eqref{key.family} and Definition \ref{Def4.3}. Given $\lambda \in (0,1]$ we write $\mathcal{F}(\lambda)$ and $\mathcal{DF}(\lambda)$
instead of $\mathcal{F}_{S}(d,\lambda)$ and $\mathcal{DF}_{S}(d,\lambda)$ respectively.
Furthermore, given $\lambda \in (0,1]$ for each $x,y \in \mathbb{R}^{n}$ we define the family
\begin{equation}
\notag
\mathcal{Q}_{x,y}(\lambda):=\{Q \ni x,y : Q \in \mathcal{F}(\lambda)\}.
\end{equation}
Now, for each $\lambda \in (0,1]$ and any $x,y \in \mathbb{R}^{n}$, we set
\begin{equation}
\label{eqq.5.1}
\widetilde{\rho}_{S,d,\lambda}(x,y):=\widetilde{\rho}_{\lambda}(x,y):=
\begin{cases}
\inf\{l(Q): Q \in \mathcal{Q}_{x,y}(\lambda)\}, \quad x \neq y \hbox{ and } \mathcal{Q}_{x,y}(\lambda) \neq \emptyset;\\
+\infty, \quad \mathcal{Q}_{x,y}(\lambda) = \emptyset;\\
0, \quad x=y.
\end{cases}
\end{equation}
For any two points $x,y \in \mathbb{R}^{n}$ such that at least one of them belongs to the set $S$ we put
\begin{equation}
\label{eq7.5}
\rho_{\lambda}(x,y):=\rho_{S,d,\lambda}(x,y):=
\inf \sum\limits_{i=0}^{N-1}\widetilde{\rho}_{\lambda}(x^{i},x^{i+1}),
\end{equation}
where the infimum is taken over all finite sets $\{x^{i}\}_{i=0}^{N} \subset \mathbb{R}^{n}$ such that
$x^{0}=x$ and $x^{N}=y$. Finally, in the case when $x,y \in \mathbb{R}^{n} \setminus S$ we define
\begin{equation}
\label{eq7.6}
\rho_{\lambda}(x,y):=\rho_{S,d,\lambda}(x,y):=
\begin{cases}
+\infty \quad \hbox{if} \quad \max\{\rho_{\lambda}(x,\xi),\rho_{\lambda}(y,\xi)\}=+\infty \quad \hbox{for some} \quad \xi \in S, \\
\max\{\|x-y\|_{\infty},\sup\limits_{\xi \in S}|\rho_{\lambda}(x,\xi)-\rho_{\lambda}(y,\xi)|\} \quad \hbox{otherwise}.
\end{cases}
\end{equation}

Recall that a pseudometric on $\mathbb{R}^{n}$
is a symmetric nonnegative function $\rho: \mathbb{R}^{n} \times \mathbb{R}^{n} \to [0,+\infty]$
satisfying the triangle inequality. Given a pseudometric $\rho$ on $\mathbb{R}^{n}$,
we will use the symbol $\mathbb{R}^{n}_{\rho}$ to denote the pseudometric space $(\mathbb{R}^{n},\rho)$. By $B^{n,\rho}_{r}(x)$
we will denote the closed ball centered in $x \in \mathbb{R}^{n}$ with radius $r$ (in pseudometric $\rho$), i.e.,
$B^{n,\rho}_{r}(x):=\{y \in \mathbb{R}^{n}:\rho(x,y) \le r\}$.

\begin{Prop}
\label{Prop7.1}
For each $\lambda \in (0,1]$ the function $\rho_{\lambda}: \mathbb{R}^{n} \times \mathbb{R}^{n} \to [0,+\infty)$ is a pseudometric on $\mathbb{R}^{n}$.
\end{Prop}

\begin{proof}
The symmetry is obvious by \eqref{eq7.5} and \eqref{eq7.6}. Furthermore, note that
\begin{equation}
\notag
\rho_{\lambda}(x,y) \geq \|x-y\|_{\infty} \quad \hbox{for any} \quad x,y \in \mathbb{R}^{n}.
\end{equation}
Hence, $\rho_{\lambda}(x,y)=0$ implies $x=y$.

It remains to verify the triangle inequality.
We fix an arbitrary triple of points $x,y,z \in \mathbb{R}^{n}$.
\textit{In the case} when $\rho_{\lambda}(x,y)=+\infty$ or $\rho_{\lambda}(y,z)=+\infty$ the triangle inequality is obvious.
Consider \textit{the case} when $\rho_{\lambda}(x,y)<+\infty$ and $\rho_{\lambda}(y,z)<+\infty$. We should consider \textit{two subcases.}
In \textit{the first subcase} at least one of 3 points $x,y,z$ (say $y$) belongs to the set $S$.
Given $\delta > 0$, let $\{x^{i}\}_{i=0}^{N}, \{x^{i}\}_{i=N+1}^{L} \subset \mathbb{R}^{n}$ be finite sets of points
such that $x^{0}=x$, $x^{N}=y$, $x^{L}=z$ and
\begin{equation}
\notag
\sum\limits_{i=0}^{N-1}\widetilde{\rho}_{\lambda}(x^{i},x^{i+1}) \le \rho_{\lambda}(x,y)+\frac{\delta}{2},
\quad \sum\limits_{i=N}^{L-1}\widetilde{\rho}_{\lambda}(x^{i},x^{i+1}) \le \rho_{\lambda}(y,z)+\frac{\delta}{2}.
\end{equation}
Summing both inequalities and using \eqref{eq7.5} we obtain
\begin{equation}
\notag
\rho_{\lambda}(x,z) \le \sum\limits_{i=0}^{L-1}\widetilde{\rho}_{\lambda}(x^{i},x^{i+1}) \le
\rho_{\lambda}(x,y)+\rho_{\lambda}(y,z)+\delta.
\end{equation}
Since $\delta>0$  was chosen arbitrarily we deduce the triangle inequality for this subcase
\begin{equation}
\label{eq7.8}
\rho_{\lambda}(x,z) \le \rho_{\lambda}(x,y)+\rho_{\lambda}(y,z).
\end{equation}
Finally, consider \textit{the second subcase} when $x,y,z \in \mathbb{R}^{n} \setminus S$.
If $\|x-z\| \geq
\sup_{\xi \in S}|\rho_{\lambda}(x,\xi)-\rho_{\lambda}(z,\xi)|$, by \eqref{eq7.6}
we get
\begin{equation}
\label{eq7.9}
\rho_{\lambda}(x,z)=\|x-z\| \le \|x-y\|+\|y-z\| \le \rho_{\lambda}(x,y)+\rho_{\lambda}(y,z).
\end{equation}
If $\|x-z\| < \sup_{\xi \in S}|\rho_{\lambda,\varepsilon}(x,\xi)-\rho_{\lambda,\varepsilon}(z,\xi)|$ by \eqref{eq7.6}, we have,
for any $\xi \in S$,
\begin{equation}
\label{eq7.10'}
\begin{split}
&|\rho_{\lambda}(x,\xi)-\rho_{\lambda}(z,\xi)|\le |\rho_{\lambda}(x,\xi)-\rho_{\lambda}(y,\xi)|
+|\rho_{\lambda}(y,\xi)-\rho_{\lambda}(z,\xi)|
\le \rho_{\lambda}(x,y)+\rho_{\lambda}(y,z).
\end{split}
\end{equation}
Taking the supremum in \eqref{eq7.10'} over all $\xi \in S$ we obtain
\begin{equation}
\label{eq7.11}
\rho_{\lambda}(x,z) \le \rho_{\lambda}(x,y)+\rho_{\lambda}(y,z).
\end{equation}
Combining \eqref{eq7.9} with \eqref{eq7.11} we get the triangle inequality for the subcase when $x,y,z \in \mathbb{R}^{n} \setminus S$

As a result, we have proved that the triangle inequality holds for any triple of points $x,y,z \in \mathbb{R}^{n}$.

This completes the proof.

\end{proof}

\begin{Remark}
\label{Rem61}
Note that if $\lambda \in (0,\lambda_{S}]$, then $Q_{0,0} \in \mathcal{DF}(\lambda)$ (by assumption ($\textbf{B}$) and Remark \ref{Rem2.1}). Hence,
by \eqref{eq7.5}, \eqref{eq7.6} it is easy to see that $\rho_{\lambda}(x,y) < +\infty$ for all $x,y \in Q_{0,0}$.
\end{Remark}
\hfill$\Box$

Given a nonempty set $E \subset \mathbb{R}^{n} \setminus S$ we define the \textit{$(d,\lambda)$-thick distance}
between the sets $E \subset \mathbb{R}^{n} \setminus S$ and $S$ by the formula
\begin{equation}
\label{eq7.13}
\operatorname{D}_{\lambda}(E,S):=\operatorname{D}_{S,d,\lambda}(E,S):=\inf \{\rho_{\lambda}(x,\xi):x \in E, \xi \in S\}.
\end{equation}

The following proposition gives a more simple way for computing the $(d,\lambda)$-thick distance
from a given point $x \in \mathbb{R}^{n}\setminus S$ to the set $S$.

\begin{Prop}
\label{Prop7.2}
Let $\lambda \in (0,\lambda_{S}]$. Then the following equality
\begin{equation}
\label{eq7.15}
\operatorname{D}_{\lambda}(x,S)=\inf \{l(Q):Q \ni x \hbox{ and $Q \in \mathcal{F}(\lambda)$}\}
\end{equation}
holds for any $x \in (\frac{\lambda_{S}}{\lambda})^{\frac{1}{d}}Q_{0,0} \setminus S$.
\end{Prop}

\begin{proof}
By Remark \ref{Rem31} and assumption (\textbf{B}) we have $(\frac{\lambda_{S}}{\lambda})^{\frac{1}{d}}Q_{0,0} \in \mathcal{F}(\lambda)$. Hence,
$\mathcal{Q}_{x,y}(\lambda) \neq \emptyset$ for each
$x \in (\frac{\lambda_{S}}{\lambda})^{\frac{1}{d}}Q_{0,0} \setminus S$ and $y \in S$.
As a result, $0 < \rho_{\lambda}(x,y) < +\infty$ for such $x$ and $y$.

Now we fix an arbitrary $x \in  (\frac{\lambda_{S}}{\lambda})^{\frac{1}{d}}Q_{0,0} \setminus S$. We denote the right-hand side
of \eqref{eq7.15} by $\widetilde{\operatorname{D}}_{\lambda}(x,S)$. Our aim is to show that $\operatorname{D}_{\lambda}(x,S)=\widetilde{\operatorname{D}}_{\lambda}(x,S)$.
Given $\varepsilon > 0$, using \eqref{eq7.5}, \eqref{eq7.13}, we choose $x_{\varepsilon} \in S$ and points $\{x_{i}\}_{i=0}^{N_{\varepsilon}}$ in such a way that
$x_{0}=x$, $x_{N_{\varepsilon}}=x_{\varepsilon}$ and
\begin{equation}
\notag
\sum\limits_{i=0}^{N_{\varepsilon}-1}\widetilde{\rho}_{\lambda}(x_{i},x_{i+1}) < \operatorname{D}_{\lambda}(x,S)+\varepsilon.
\end{equation}
Hence, we clearly get
\begin{equation}
\notag
\operatorname{D}_{\lambda}(x,S) \le \widetilde{\operatorname{D}}_{\lambda}(x,S) \le \widetilde{\rho}_{\lambda}(x_{0},x_{1}) < \operatorname{D}_{\lambda}(x,S)+\varepsilon.
\end{equation}
Since $\varepsilon > 0$ can be chosen arbitrary small, we get the required equality and complete the proof.
\end{proof}

\begin{Remark}
\label{Rem7.2}
The pseudometric $\rho_{\lambda}:=\rho_{S,d,\lambda}$ introduced above is a natural generalization of the metric induced by the $\|\cdot\|_{\infty}$-norm.
Indeed, given $\lambda \in (0,1]$ we have $\rho_{S,0,\lambda}(x,y)=\|x-y\|_{\infty}$ for any $x,y \in \mathbb{R}^{n}$.

Given $\lambda \in (0,\lambda_{S}]$, recall the concept of $(d,\lambda)$-thick sets in $\mathbb{R}^{n}$ formulated in the introduction. Recall also the notion of
$\delta$-neighborhood of $S$. If $S$ is $(d,\lambda)$-thick, then using Proposition \ref{Prop7.2}
it is easy to see that
\begin{equation}
\notag
\operatorname{dist}(x,S) \le \operatorname{D}_{\lambda}(x,S) \le 2\operatorname{dist}(x,S) \quad \hbox{for all} \quad x \in U_{\frac{1}{2}}(S).
\end{equation}
\end{Remark}
\hfill$\Box$

Given $\lambda \in (0,1]$, for any $\varepsilon > 0$ we introduce the $\varepsilon$-neighborhood of the set $S$ with respect to the
metric $\rho_{\lambda}$ by the formula
\begin{equation}
\notag
U^{\rho_{\lambda}}_{\varepsilon}(S):=\{x \in \mathbb{R}^{n}: \operatorname{D}_{\lambda}(x,S) < \varepsilon\}.
\end{equation}
By Remark \ref{Remm2.4} we have
\begin{equation}
\notag
U^{\rho_{\lambda}}_{\varepsilon}(S) \supset S \quad \hbox{for all} \quad \varepsilon \geq 0 \quad \hbox{and all} \quad \lambda \in (0,1].
\end{equation}

\begin{Remark}
For each $\lambda \in (0,\lambda_{S}]$ and $\varepsilon > 0$ small enough, we have (recall \eqref{eqq.1.2})
\begin{equation}
U^{\rho_{\lambda}}_{\varepsilon}(S)=S_{\varepsilon}(\lambda):=\cup \{Q: Q \in \mathcal{F}(\lambda) \hbox{ and } 0 \le l(Q) < \varepsilon\}.
\end{equation}
Indeed, by Proposition \ref{Prop7.2} for a given point $x \in (\frac{\lambda_{S}}{\lambda})^{d}Q_{0,0} \setminus S$ we have
$\operatorname{D}_{\lambda}(x,S) < \varepsilon$ if and only if there is a cube $Q \in \mathcal{F}(\lambda)$ with $0 \le l(Q) < \varepsilon$.
This proves the claim.
\end{Remark}

\hfill$\Box$

\textit{Proof of Theorem \ref{Th.cavity.decom.}.} We set
\begin{equation}
\overline{\lambda}:=\frac{\lambda_{S}}{3^{n}c^{d}}
\end{equation}
and fix $\lambda \in (0,\overline{\lambda}]$. We split the proof into several steps.

\textit{Step 1.}
Note that by assumption (\textbf{B})  it follows that $Q_{0,0} \in \mathcal{F}(\lambda_{S})$. Hence, by Remark \ref{Rem31}
we deduce that the cube $cQ_{0,0} \in \mathcal{F}(\frac{\lambda_{S}}{c^{d}})$. By Proposition \ref{Prop7.2},
\begin{equation}
\label{eq718}
\operatorname{D}_{\lambda}(x,S) \le c \quad \hbox{for all} \quad x \in cQ_{0,0} \setminus S.
\end{equation}

For each $k \in Z$ and every $i \in \mathbb{N}_{0}$ we define the layers
\begin{equation}
\label{eq719}
\operatorname{L}_{k,i}:=\Bigl\{x \in cQ_{0,0} \setminus S: \frac{c}{2^{k+i+1}} \le \operatorname{D}_{\lambda}(x,S) < \frac{c}{2^{k}}\Bigr\}.
\end{equation}
We set $\operatorname{L}_{k}:=\operatorname{L}_{k,0}(\lambda)$ for brevity.
By \eqref{eq718} and \eqref{eq719}  we have
\begin{equation}
\label{eqq.514}
\operatorname{int}(cQ_{0,0}) \setminus S=\cup_{k \in \mathbb{N}_{0}}\operatorname{L}_{k,i}(\lambda) \quad
\hbox{for every} \quad i \in \mathbb{N}_{0}.
\end{equation}

\textit{Step 2.}
For each $x \in cQ_{0,0} \setminus S$, there
exists a unique number $k_{x} \in \mathbb{N}_{0}$
for which $x \in \operatorname{L}_{k_{x}}$.
Furthermore, by \eqref{eq719} and Proposition \ref{Prop7.2} there is a cube $Q_{x} \in \mathcal{F}(\lambda)$ such that
\begin{equation}
\label{eqq.515}
x \in Q_{x} \subset cQ_{0,0} \quad \hbox{and} \quad \frac{c}{2^{k_{x}+1}} \le l(Q_{x}) < \frac{c}{2^{k_{x}}}.
\end{equation}
By \eqref{eqq.515} and Proposition \ref{Prop7.2} it is clear that
\begin{equation}
\label{eqq.516}
\operatorname{D}_{\lambda}(y,S) < \frac{c}{2^{k_{x}}} \quad \hbox{for all} \quad y \in Q_{x}.
\end{equation}

For each point $x \in cQ_{0,0} \setminus S$, we find (and fix to the end of the proof)
an arbitrary pair $\alpha_{x}:=(k_{x}+4,m_{x}) \in \mathbb{N} \times \mathbb{Z}^{n}$ such that
\begin{equation}
\label{eqq.5.17}
Q_{\alpha_{x}} \subset Q_{x} \quad \hbox{and} \quad x \in 3Q_{\alpha_{x}}.
\end{equation}
We have
\begin{equation}
\label{eqq.5.20}
3Q_{\alpha_{x}} \notin \mathcal{F}(\lambda).
\end{equation}
Indeed, otherwise
by Proposition \ref{Prop7.2} we would get
\begin{equation}
\notag
\operatorname{D}_{\lambda}(x,S) \le l(3Q_{\alpha_{x}}) < \frac{1}{4}2^{-k_{x}} < \frac{c}{2^{k_{x}+1}}.
\end{equation}
This inequality implies that $x \notin \operatorname{L}_{k_{x}}$ contradicting the definition of $k_{x}$.
From \eqref{eqq.5.20} and Remark \ref{Rem31} it follows that
\begin{equation}
\label{eqq.5.21}
Q_{\alpha_{x}} \notin \mathcal{F}(\frac{\lambda_{S}}{c^{d}}).
\end{equation}

\textit{Step 3.}
Let $\overline{\delta}:=\overline{\delta}(n,d,\lambda,\frac{\lambda_{S}}{c^{d}}) \in (0,1)$ be the same as in Theorem \ref{Th1.1}. We set
\begin{equation}
\notag
j^{\ast}:=j^{\ast}(n,d,\lambda):=-[\log_{2}\frac{\overline{\delta}}{c}]+5 \in \mathbb{N}.
\end{equation}
For each $x \in cQ_{0,0} \setminus $ we put
\begin{equation}
\label{eq722}
\Omega_{\alpha_{x}}:= 3Q_{\alpha_{x}} \cap  \operatorname{L}_{k_{x},j^{\ast}+1}.
\end{equation}
Finally, we define the family
\begin{equation}
\notag
\mathcal{W}_{S}:=\mathcal{W}_{S}(d,\lambda,c):=\{\Omega_{\alpha_{x}}: x \in cQ_{0,0} \setminus S\}.
\end{equation}

\textit{Step 4.} Now we show that the family $\mathcal{W}_{S}$ has all necessary properties.
Since $x \in 3Q_{\alpha_{x}} \cap \operatorname{L}_{k_{x}}$ we clearly have
\begin{equation}
x \in \Omega_{\alpha_{x}}.
\end{equation}
Hence,
\begin{equation}
\label{eq725}
\operatorname{int}(cQ_{0,0}) \setminus S = \cup \{\Omega: \Omega \in \mathcal{W}_{S}\}.
\end{equation}
This proves assertion (i) of the theorem.

It is clear that
\begin{equation}
\label{eq726}
M(\{\operatorname{L}_{k,i}(\lambda)\}_{k \in \mathbb{N}_{0}}) \le i+1 \quad \hbox{for every} \quad i \in \mathbb{N}_{0}.
\end{equation}
On the other hand, $\Omega_{\alpha_{x}} \subset 3Q_{\alpha_{x}}$ for all $x \in cQ_{0,0} \setminus S$.
Combining this fact with Proposition \ref{Prop21} and \eqref{eq726} we deduce that
\begin{equation}
\label{eq727}
M(\mathcal{W}_{S}) \le M(\{\operatorname{L}_{k,j^{\ast}+1}(\lambda)\}_{k \in \mathbb{N}_{0}})M(3\mathcal{D}_{k}) \le 4^{n}(j^{\ast}+1).
\end{equation}
This proves assertion (iv) of the theorem.

Note that by the construction we have
\begin{equation}
\label{eq525'''}
\frac{c}{2^{k_{x}+j^{\ast}}} \le \operatorname{D}_{\lambda}(y,S) < \frac{c}{2^{k_{x}}} \quad \hbox{for all} \quad y \in \Omega_{\alpha_{x}}.
\end{equation}
Combining this fact with \eqref{eqq.516} and taking into account that $Q_{\alpha_{x}} \subset Q_{x}$, we obtain
\begin{equation}
\label{eq.forgot}
Q_{\alpha_{x}} \setminus \bigcup_{j=j^{\ast}+1}^{\infty}\operatorname{L}_{k_{x}+j} \subset \Omega_{\alpha_{x}}.
\end{equation}
Keeping in mind that $2^{-k_{x}-j^{\ast}} < \overline{\delta}l(Q_{\alpha_{x}})$ we collect
\eqref{eq719}, \eqref{eqq.5.21}, \eqref{eq.forgot}, \eqref{eq722} and apply Theorem \ref{Th1.1} with $\overline{Q}=Q_{\alpha_{x}}$ we deduce
\begin{equation}
\label{eq729}
\mathcal{L}^{n}(\Omega_{\alpha_{x}}) \geq  \underline{\gamma}(l(Q_{\alpha_{x}}))^{n} \geq \frac{\underline{\gamma}}{3^{n}}(\operatorname{diam}\Omega_{\alpha_{x}})^{n},
\end{equation}
where  $\underline{\gamma}=\underline{\gamma}(n,d,\frac{\lambda_{S}}{c^{d}})$ is the same as in Theorem \ref{Th1.1}.
This proves assertion (iii) of the theorem.

Furthermore, from the first inequality in \eqref{eq729} it follows immediately that
\begin{equation}
\label{eq729'''}
\operatorname{diam}\Omega_{\alpha_{x}} \geq (\underline{\gamma})^{\frac{1}{n}}l(Q_{\alpha_{x}}).
\end{equation}
Combining \eqref{eq525'''}, \eqref{eq729}, \eqref{eq729'''} we prove assertion (ii) of the theorem.

The proof is complete.
\hfill$\Box$

\begin{Remark}
Analysis of the proof of Theorem \ref{Th.cavity.decom.} shows that in fact
for every $\Omega \in \mathcal{W}_{S}$,
\begin{equation}
\frac{1}{\operatorname{C}_{1}}\operatorname{diam}\Omega
\le \operatorname{D}_{S,d,\lambda}(y,S) \le \operatorname{C}_{1}\operatorname{diam}\Omega \quad \hbox{for all} \quad y \in \Omega.
\end{equation}
\end{Remark}
\hfill$\Box$




\section{Examples}

In the concluding section of the present paper we present elementary examples demonstrating the sharpness of
the main results.

The following example shows
that the $d$-thick condition in Theorem \ref{Ca.1.1} cannot be dropped.

\begin{Example}
\label{ex6.1}
Let $K \subset [0,1]$ be the standard middle-third Cantor set. For each $j \in \mathbb{N}$ we define $K_{j}:=\{3^{-j}x: x \in K\}$. We set
\begin{equation}
\notag
S_{j}:=\bigcup_{i=0}^{2^{j}-1}\Bigl(\frac{i}{2^{j}}+K_{2j}\Bigr).
\end{equation}
Obviously, given $j \in \mathbb{N}$, the maximal size of closed $1$-dimensional cubes
$Q \subset [0,1] \setminus S_{j}$ is less than $2^{-j}$.
If $d=\operatorname{ln}2/\operatorname{ln}3$ we clearly have
\begin{equation}
\label{eq6.1}
\mathcal{H}^{d}(S_{j}) \le \frac{2^{j}}{2^{2j}} \to 0, \quad j \to \infty.
\end{equation}
Finally, we define the set
\begin{equation}
\notag
S:= \bigcup_{j=0}^{\infty} \Bigl((1-\frac{1}{2^{j}})+\frac{1}{2^{j+1}}S_{j+1}\Bigr).
\end{equation}
It follows from \eqref{eq6.1} that the set $S$ is not $(d,\lambda)$-thick for any $\lambda \in (0,1)$.
On the other hand, any dyadic interval $[1-\frac{2}{2^{j}},1-\frac{1}{2^{j}}]$, $j \in \mathbb{N}$, can be $(S,\tau)$-porous only with $\tau < 2^{-s}$.
\end{Example}

\hfill$\Box$

Now we show that the restriction $d < n$  in Theorem \ref{Th1.1} can not be dropped.

\begin{Example}
\label{ex6.3}
For each $j \in \mathbb{N}$ we set
\begin{equation}
S_{j}:=\bigcup_{i=0}^{2^{j}-1}[\frac{i}{2^{j}},\frac{i}{2^{j}}+\frac{1}{10}\frac{1}{2^{j}}].
\end{equation}
Finally, we define the set
\begin{equation}
\notag
S:= \bigcup_{j=0}^{\infty} \Bigl((1-\frac{1}{2^{j}})+\frac{1}{2^{j+1}}S_{j+1}\Bigr).
\end{equation}
It is easy to see that
\begin{equation}
\mathcal{H}^{1}_{\infty}([1-\frac{2}{2^{j}},1-\frac{1}{2^{j}}] \cap S) < \frac{1}{8}\frac{1}{2^{j}} \quad \hbox{for all} \quad j \in \mathbb{N}.
\end{equation}
On the other hand, for any $j \in \mathbb{N}$,
\begin{equation}
\Bigl[1-\frac{2}{2^{j}},1-\frac{1}{2^{j}}\Bigr] \setminus \cup\{Q':l(Q') \le 2^{-2j} \hbox{ and $Q' \in \mathcal{F}_{S}(1,\frac{1}{10})$}\} = \emptyset.
\end{equation}
This shows that the conclusion of Theorem \ref{Th1.1} cannot hold in the case $d=n=1$ and $\lambda = \frac{1}{10}$.
\end{Example}

\hfill$\Box$

\end{document}